\documentclass[12pt]{article}
\usepackage{latexsym, amssymb}
\DeclareMathAlphabet\gothic{U}{euf}{m}{n}

\makeatletter
\def\eqnarray{\stepcounter{equation}\let\@currentlabel=\theequation
\global\@eqnswtrue
\tabskip\@centering\let\\=\@eqncr
$$\halign to \displaywidth\bgroup\hfil\global\@eqcnt\z@
  $\displaystyle\tabskip\z@{##}$&\global\@eqcnt\@ne
  \hfil$\displaystyle{{}##{}}$\hfil
  &\global\@eqcnt\tw@ $\displaystyle{##}$\hfil
  \tabskip\@centering&\llap{##}\tabskip\z@\cr}

\def\endeqnarray{\@@eqncr\egroup
      \global\advance\c@equation\m@ne$$\global\@ignoretrue}

\def\@yeqncr{\@ifnextchar [{\@xeqncr}{\@xeqncr[5pt]}}
\makeatother

\begin{document}
\bibliographystyle{tom}

\newtheorem{lemma}{Lemma}[section]
\newtheorem{thm}[lemma]{Theorem}
\newtheorem{cor}[lemma]{Corollary}
\newtheorem{voorb}[lemma]{Example}
\newtheorem{rem}[lemma]{Remark}
\newtheorem{prop}[lemma]{Proposition}
\newtheorem{stat}[lemma]{{\hspace{-5pt}}}
\newtheorem{obs}[lemma]{Observation}
\newtheorem{defin}[lemma]{Definition}

\newenvironment{remarkn}{\begin{rem} \rm}{\end{rem}}
\newenvironment{exam}{\begin{voorb} \rm}{\end{voorb}}
\newenvironment{defn}{\begin{defin} \rm}{\end{defin}}
\newenvironment{obsn}{\begin{obs} \rm}{\end{obs}}

\newenvironment{emphit}{\begin{itemize} }{\end{itemize}}

\newcommand{\gota}{\gothic{a}}
\newcommand{\gotb}{\gothic{b}}
\newcommand{\gotc}{\gothic{c}}
\newcommand{\gote}{\gothic{e}}
\newcommand{\gotf}{\gothic{f}}
\newcommand{\gotg}{\gothic{g}}
\newcommand{\gothh}{\gothic{h}}
\newcommand{\gotk}{\gothic{k}}
\newcommand{\gotm}{\gothic{m}}
\newcommand{\gotn}{\gothic{n}}
\newcommand{\gotp}{\gothic{p}}
\newcommand{\gotq}{\gothic{q}}
\newcommand{\gotr}{\gothic{r}}
\newcommand{\gots}{\gothic{s}}
\newcommand{\gotu}{\gothic{u}}
\newcommand{\gotv}{\gothic{v}}
\newcommand{\gotw}{\gothic{w}}
\newcommand{\gotz}{\gothic{z}}
\newcommand{\gotA}{\gothic{A}}
\newcommand{\gotB}{\gothic{B}}
\newcommand{\gotG}{\gothic{G}}
\newcommand{\gotL}{\gothic{L}}
\newcommand{\gotS}{\gothic{S}}
\newcommand{\gotT}{\gothic{T}}

\newcommand{\mn}{\marginpar{\hspace{1cm}*} }
\newcommand{\mnn}{\marginpar{\hspace{1cm}**} }

\newcommand{\mnq}{\marginpar{\hspace{1cm}*???} }
\newcommand{\mnnq}{\marginpar{\hspace{1cm}**???} }

\newcounter{teller}
\renewcommand{\theteller}{\Roman{teller}}
\newenvironment{tabel}{\begin{list}%
{\rm \bf \Roman{teller}.\hfill}{\usecounter{teller} \leftmargin=1.1cm
\labelwidth=1.1cm \labelsep=0cm \parsep=0cm}
                      }{\end{list}}

\newcounter{tellerr}
\renewcommand{\thetellerr}{(\roman{tellerr})}
\newenvironment{subtabel}{\begin{list}%
{\rm  (\roman{tellerr})\hfill}{\usecounter{tellerr} \leftmargin=1.1cm
\labelwidth=1.1cm \labelsep=0cm \parsep=0cm}
                         }{\end{list}}
\newenvironment{ssubtabel}{\begin{list}%
{\rm  (\roman{tellerr})\hfill}{\usecounter{tellerr} \leftmargin=1.1cm
\labelwidth=1.1cm \labelsep=0cm \parsep=0cm \topsep=1.5mm}
                         }{\end{list}}

\newcommand{\Ni}{{\bf N}}
\newcommand{\Ri}{{\bf R}}
\newcommand{\Ci}{{\bf C}}
\newcommand{\Qi}{{\bf Q}}
\newcommand{\Ti}{{\bf T}}
\newcommand{\Zi}{{\bf Z}}
\newcommand{\Fi}{{\bf F}}

\newcommand{\proof}{\mbox{\bf Proof} \hspace{5pt}} 
\newcommand{\remark}{\mbox{\bf Remark} \hspace{5pt}}
\newcommand{\ruimte}{\vskip10.0pt plus 4.0pt minus 6.0pt}

\newcommand{\simh}{{\stackrel{{\rm cap}}{\sim}}}
\newcommand{\ad}{{\mathop{\rm ad}}}
\newcommand{\Ad}{{\mathop{\rm Ad}}}
\newcommand{\Aut}{\mathop{\rm Aut}}
\newcommand{\arccot}{\mathop{\rm arccot}}
\newcommand{\capp}{{\mathop{\rm cap}}}
\newcommand{\Capp}{{\mathop{\rm Cap}}}
\newcommand{\rcapp}{{\mathop{\rm rcap}}}
\newcommand{\diam}{\mathop{\rm diam}}
\newcommand{\divv}{\mathop{\rm div}}
\newcommand{\codim}{\mathop{\rm codim}}
\newcommand{\RRe}{\mathop{\rm Re}}
\newcommand{\IIm}{\mathop{\rm Im}}
\newcommand{\Tr}{{\mathop{\rm Tr}}}
\newcommand{\Vol}{{\mathop{\rm Vol}}}
\newcommand{\card}{{\mathop{\rm card}}}
\newcommand{\supp}{\mathop{\rm supp}}
\newcommand{\sgn}{\mathop{\rm sgn}}
\newcommand{\essinf}{\mathop{\rm ess\,inf}}
\newcommand{\esssup}{\mathop{\rm ess\,sup}}
\newcommand{\Int}{\mathop{\rm Int}}
\newcommand{\Leibniz}{\mathop{\rm Leibniz}}
\newcommand{\lcm}{\mathop{\rm lcm}}
\newcommand{\loc}{{\rm loc}}

\newcommand{\mod}{\mathop{\rm mod}}
\newcommand{\spann}{\mathop{\rm span}}
\newcommand{\one}{1\hspace{-4.5pt}1}

\newcommand{\DWR}{}

\hyphenation{groups}
\hyphenation{unitary}

\newcommand{\tfrac}[2]{{\textstyle \frac{#1}{#2}}}

\newcommand{\cb}{{\cal B}}
\newcommand{\cc}{{\cal C}}
\newcommand{\cd}{{\cal D}}
\newcommand{\ce}{{\cal E}}
\newcommand{\cf}{{\cal F}}
\newcommand{\ch}{{\cal H}}
\newcommand{\ci}{{\cal I}}
\newcommand{\ck}{{\cal K}}
\newcommand{\cl}{{\cal L}}
\newcommand{\cm}{{\cal M}}
\newcommand{\cn}{{\cal N}}
\newcommand{\co}{{\cal O}}
\newcommand{\cs}{{\cal S}}
\newcommand{\ct}{{\cal T}}
\newcommand{\cx}{{\cal X}}
\newcommand{\cy}{{\cal Y}}
\newcommand{\cz}{{\cal Z}}

\newcommand{\wtozp}{W^{1,2}\raisebox{10pt}[0pt][0pt]{\makebox[0pt]{\hspace{-34pt}$\scriptstyle\circ$}}}
\newlength{\hightcharacter}
\newlength{\widthcharacter}
\newcommand{\covsup}[1]{\settowidth{\widthcharacter}{$#1$}\addtolength{\widthcharacter}{-0.15em}\settoheight{\hightcharacter}{$#1$}\addtolength{\hightcharacter}{0.1ex}#1\raisebox{\hightcharacter}[0pt][0pt]{\makebox[0pt]{\hspace{-\widthcharacter}$\scriptstyle\circ$}}}
\newcommand{\cov}[1]{\settowidth{\widthcharacter}{$#1$}\addtolength{\widthcharacter}{-0.15em}\settoheight{\hightcharacter}{$#1$}\addtolength{\hightcharacter}{0.1ex}#1\raisebox{\hightcharacter}{\makebox[0pt]{\hspace{-\widthcharacter}$\scriptstyle\circ$}}}
\newcommand{\scov}[1]{\settowidth{\widthcharacter}{$#1$}\addtolength{\widthcharacter}{-0.15em}\settoheight{\hightcharacter}{$#1$}\addtolength{\hightcharacter}{0.1ex}#1\raisebox{0.7\hightcharacter}{\makebox[0pt]{\hspace{-\widthcharacter}$\scriptstyle\circ$}}}

 \thispagestyle{empty}
 
 \begin{center}

{\Large{\bf $L_1$-uniqueness }}\\[3mm] 
{\Large{\bf of degenerate elliptic operators  }}  \\[5mm]
\large  Derek W. Robinson$^1$ and Adam Sikora$^2$\\[2mm]

\normalsize{September  2010}
\end{center}

\vspace{5mm}

\begin{center}
{\bf Abstract}
\end{center}

\begin{list}{}{\leftmargin=1.8cm \rightmargin=1.8cm \listparindent=15mm 
   \parsep=0pt}
   \item
   Let $\Omega$ be an open subset of $\Ri^d$ with $0\in \Omega$.
   Further let  $H_\Omega=-\sum^d_{i,j=1}\partial_i\,c_{ij}\,\partial_j$
be a   second-order partial differential operator  with domain $C_c^\infty(\Omega)$  
where  the coefficients $c_{ij}\in W^{1,\infty}_{\rm loc}(\overline\Omega)$ are real, $c_{ij}=c_{ji}$ and  the
coefficient  matrix $C=(c_{ij})$ satisfies bounds $0<C(x)\leq c(|x|) I$ for all $x\in \Omega$.
If \[
\int^\infty_0ds\,s^{d/2}\,e^{-\lambda\,\mu(s)^2}<\infty
\]
for some  $\lambda>0$ where $\mu(s)=\int^s_0dt\,c(t)^{-1/2}$ then we establish
 that  $H_\Omega$ is $L_1$-unique, i.e.\ it has a unique
$L_1$-extension   which generates a continuous semigroup, if and only if it is Markov unique,  
i.e.\ it has a unique $L_2$-extension   which generates a submarkovian semigroup.
Moreover these uniqueness conditions are equivalent with the capacity of the boundary of $\Omega$,
measured with respect to $H_\Omega$,  being zero.
We also demonstrate that the capacity  depends on two gross features,
the Hausdorff dimension of subsets $A$ of the boundary the set and  the order of degeneracy of $H_\Omega$ at~$A$.

\end{list}

\vfill


\vspace{0.5cm}

\noindent
\begin{tabular}{@{}cl@{\hspace{10mm}}cl}
1. & Centre for Mathematics & 
  2. & Department of Mathematics  \\
&\hspace{15mm} and its Applications  & 
  &Macquarie University  \\
& Mathematical Sciences Institute & 
  & Sydney, NSW 2109  \\
& Australian National University& 
  & Australia \\
& Canberra, ACT 0200  & {}
  & \\
& Australia & {}
  & \\
& derek.robinson@anu.edu.au & {}
  &sikora@ics.mq.edu.au \\
\end{tabular}

\newpage

\setcounter{page}{1}

\section{Introduction}\label{S1}

In a recent paper  \cite{RSi4} we established that Markov uniqueness and $L_1$-uniqueness
are equivalent properties  for  a second-order, symmetric, elliptic operator with bounded Lipschitz continuous coefficients
$c_{ij}$ on an open subset $\Omega$ of $\Ri^d$.
Moreover, these  properties 
hold if and only if the corresponding capacity of the boundary 
$\partial\Omega$ of $\Omega$ is zero.
In this note we extend these results to operators with locally bounded coefficients with a possible growth at infinity. 
As an illustration of our results  we establish that    Markov uniqueness, $L_1$-uniqueness and the capacity condition
are equivalent if the matrix $C=(c_{ij})$ satisfies  $\|C(x)\|\sim |x|^2(\log |x|)^\alpha$ as $|x|\to\infty$
with $\alpha\in[0,1]$.
In addition we give an example with  $\|C(x)\|\sim |x|^2(\log |x|)^{1+ \varepsilon}$ as $|x|\to\infty$,
where $\varepsilon>0$ is  arbitrarily small,  which is  Markov unique but not $L_1$-unique.
Our results extend uniqueness criteria previously established for the special case $\Omega=\Ri^d$
(see \cite{Dav14}, \cite{Ebe} Chapter~2, \cite{Sta} Section~2, and references therein).

Let $\Omega$ be an open subset of $\Ri^d$ and choose coordinates such that $0\in\Omega$.  
Define $H_\Omega$ as  the positive symmetric operator  on $L_2(\Omega)$ with
 domain $D(H_\Omega)=C_c^\infty(\Omega)$ and action
\begin{equation}
H_\Omega\varphi= -\sum^d_{i,j=1}\partial_i\,c_{ij}\,\partial_j\varphi
\label{eubc1.0}
\end{equation}
where  $\partial_i=\partial/\partial x_i$ and 
 the coefficients $c_{ij}$ satisfy
\begin{equation}
\left.
\begin{array}{ll}
1. & c_{ij}=c_{ji}\in W^{1,\infty}_{\rm loc}(\overline \Omega) \mbox{ are real, }\\[10pt]
2. &C(x)=(c_{ij}(x))
\mbox{ is a  strictly positive-definite matrix for all } x\in \Omega
 \end{array}
 \right\}
 \label{eubc1.1}
 \end{equation}
 where $W^{s,p}_{\rm loc}(\Omega)$ denotes the local version of the usual Sobolev spaces
 and $W^{s,p}_{\rm loc}(\overline\Omega)$ denotes the restriction to $\Omega$ of functions in 
 $W^{s,p}_{\rm loc}(\Ri^d)$. 
 The class of operators defined by (\ref{eubc1.0}) and (\ref{eubc1.1}) will be denoted by $\ce_\Omega$.
 
 It follows that each $H_\Omega\in \ce_\Omega$ is locally strongly elliptic, i.e.\
 for each relatively compact $V\subset\Omega$ there are $\mu_V,\lambda_V>0$ such that 
 $\mu_V I\leq C(x)\leq \lambda_V I$ for all $x\in V$.
There are, however, two potential sources of degeneracy.
 It is possible that $c_{ij}(x)\to0$ as $x\to \partial\Omega$ or that $c_{ij}(x)\to\infty$ as $|x|\to\infty$.

In order to control the possible growth of the coefficients at infinity we introduce
 the strictly positive non-decreasing function $c$ by
\begin{equation}
 r\in\langle0,\infty\rangle\mapsto c(r)=\sup\{\|C(x)\|: x\in \Omega , |x|<r\}
 \label{eubc1.10}
 \end{equation}
 where $\|C(x)\|$ denotes the norm of the matrix $C(x)=(c_{ij}(x))$.
 Then $\|C(x)\|\leq c(|x|)$ and $c(0_+)>0$. 
The growth conditions will be expressed either explicitly or implicitly  in terms of the asymptotic properties of the positive increasing function $\mu$ given~by
\begin{equation}
 s\in\langle0,\infty\rangle\mapsto\mu(s)=\int^s_0dt\, c(t)^{-1/2}
 \;.
 \label{eubc1.11}
 \end{equation}
 This function  is  a lower bound on the  Riemannian  distance to infinity measured with respect to   the metric $C^{-1}$.
If, for example,  $c(s)\sim s^2\,(\log s)^\alpha$ as $s\to\infty$ with $\alpha\in[0,2\rangle$
then $\mu(s)\sim (\log s)^{1-\alpha/2}\to\infty$ as $s\to\infty$.

We are interested in criteria for various uniqueness properties of $H_\Omega$ and adopt the terminology of \cite{Ebe}.
In particular $H_\Omega$, viewed as an operator on $L_p(\Omega)$ for  $p\in[1,\infty]$, is defined to be $L_p$-unique if it has a unique extension which generates an $L_p$-continuous semigroup.
Moreover, it is defined to be Markov unique if it has a unique self-adjoint extension on $L_2(\Omega)$ which generates a submarkovian semigroup,
i.e.\  an $L_2$-continuous contraction semigroup $S$ with the property that $0\leq S_t\varphi\leq \one$ whenever  $0\leq \varphi\leq \one$.
It follows that $H_\Omega$ is $L_2$-unique if and only if it is essentially self-adjoint (see \cite{Ebe}, Corollary~1.1.2). 
Then the self-adjoint closure is automatically submarkovian and  $H_\Omega$ is Markov unique.
Moreover, if $H_\Omega$ is $L_1$-unique then it is Markov unique (\cite{Ebe},  Lemma~1.1.6).

First, introduce the positive quadratic form $h_\Omega$  associated with $H_\Omega$ by 
\[
h_\Omega(\varphi)
=\sum^d_{i,j=1} \int_\Omega dx\, c_{ij}(x)\,(\partial_i\varphi)(x)(\partial_j\varphi)(x)=(\varphi,H_\Omega\varphi)
\]
with domain $D(h_\Omega)=D(H_\Omega)=C_c^\infty(\Omega)$.
Since $h_\Omega$ is the form of the symmetric operator $H_\Omega$ it  is closable with respect to the graph norm $\varphi\mapsto\|\varphi\|_{D(h_\Omega)}=(h_\Omega(\varphi)+\|\varphi\|_2^2)^{1/2}$.
In the sequel we  use the well known relationship between positive closed quadratic forms and positive self-adjoint operators  (see \cite{Kat1}, Chapter~6) together with  the corresponding theory
of Dirichlet forms and submarkovian operators (see \cite{BH} \cite{MR} \cite{FOT}).
The closure $\overline h_\Omega$ of $h_\Omega$ is automatically a Dirichlet form and the corresponding
 positive self-adjoint operator,  the Friedrichs extension $H_\Omega^F$  of $H_\Omega$, 
is submarkovian. 
Formally  $H_\Omega^F$  corresponds to the self-adjoint extension of $H_\Omega$ with Dirichlet  conditions on the boundary $\partial\Omega$ of $\Omega$.
In order to emphasize this interpretation we  adopt the alternative notation $H_{\Omega,D}= H_\Omega^F$ and $h_{\Omega, D}=\overline h_\Omega$. 

Secondly,  we  introduce  a positive self-adjoint extension of $H_\Omega$ related to Neumann boundary conditions.
Let  $\chi\in C_c^\infty(\Omega)$ with $0\leq \chi\leq \one_\Omega$  and define 
$h_{\Omega,\chi}$ as  the form of the symmetric operator on $L_2(\Omega)$ with coefficients $\chi \,c_{ij}$.
Then $h_{\Omega,\chi}$ is closable,  its closure $\overline h_{\Omega,\chi}$ is a Dirichlet form and  $\overline h_{\Omega,\chi}\leq h_{\Omega, D}$.
Next set $\cc_\Omega=\{\chi\in C_c^\infty(\Omega), \; 0\leq \chi\leq \one_\Omega\}$.
It follows that $\cc_\Omega$  is a convex  set  which is directed with respect to the natural order
and  if $\chi, \eta\in \cc_\Omega$ with $\chi\leq \eta$  then $ \overline h_{\Omega,\chi}\leq \overline h_{\Omega,\eta}$.
Now we define $h_{\Omega,N}$ by 
\begin{equation}
h_{\Omega,N}(\varphi)=\lim\{\overline h_{\Omega,\chi}(\varphi): \chi\in \cc_\Omega\}
=\sup\{\overline h_{\Omega,\chi}(\varphi): \chi\in \cc_\Omega\}
\;.
\label{emin1.1}
\end{equation}
Since $h_{\Omega,N}$ is the limit of quadratic forms it is a quadratic form
and since it is the supremum of a family of closed forms it is a closed form.
It is automatically a Dirichlet form satisfying $h_{\Omega,N}\leq h_{\Omega, D}$.
If  $H_{\Omega,N}$ is the positive self-adjoint operator associated
with $h_{\Omega,N}$ it readily follows that  $H_{\Omega,N}$ is a submarkovian extension of $H_\Omega$
and $H_{\Omega, N}\leq H_{\Omega, D}$.
If $\partial\Omega$ is smooth, or even Lipschitz continuous, then $H_{\Omega,N}$ corresponds to the extension of $H_\Omega$ with Neumann boundary conditions
but we adopt this definition for general open $\Omega$.

In order to formulate our main result on uniqueness properties we need two extra definitions.

The operator 
$H_\Omega\in\ce_\Omega$ is defined to be  conservative if 
the submarkovian semigroup $S^{\Omega, D}$ generated by $H_{\Omega,D}$ is conservative, 
i.e.\  if $S^{\Omega, D}_t\one_\Omega=\one_\Omega$ for all $t>0$.
Moreover, the  capacity   of the measurable subset $A\subset \overline \Omega$ relative to the operator $H_\Omega$  is defined by
\begin{eqnarray*}
\capp_\Omega(A)=\inf\Big\{\;\|\psi\|_{D(h_{\Omega,N})}^2&&\;: \;\psi\in D(h_{\Omega,N}) \mbox{ and  there exists   an open set  }\nonumber\\[-5pt]
&& U\subset \Ri^d
\mbox{ such that } U\supseteq A
\mbox{ and } \psi\geq1 \mbox{ a.\ e.\ on } U\cap\Omega\;\Big\}
\;.
\end{eqnarray*}
Thus $\capp_\Omega$ corresponds to the capacity relative to the Dirichlet form $h_{\Omega,N}$ 
as defined in \cite{BH}  or \cite{FOT}.

\begin{thm}\label{tsm1.1}
Assume  $H_\Omega\in \ce_\Omega$.
Consider the 
 following conditions:
\begin{tabel}
\item\label{tsm1.1-1}
$H_\Omega$ is conservative,
\item\label{tsm1.1-2}
$H_\Omega$ is $L_1$-unique, 
\item\label{tsm1.1-3}
$H_\Omega$ is Markov unique, 
\item\label{tsm1.1-4}
$\capp_\Omega(\partial\Omega)=0$.
 \end{tabel}
 
 Then {\rm  \ref{tsm1.1-1}$\Leftrightarrow$\ref{tsm1.1-2}$\Rightarrow$\ref{tsm1.1-3}$\Rightarrow$\ref{tsm1.1-4}}.

\medskip

Conversely, if  $\mu(s)\to\infty$ as $s\to\infty$
then~{\rm \ref{tsm1.1-4}$\Rightarrow$\ref{tsm1.1-3}}
or   if 
\begin{equation}
\int^\infty_0ds\, s^{d/2}\,e^{-\lambda \,\mu(s)^2}<\infty
\label{eubc1.2}
\end{equation}
for one   $\lambda>0$ then {\rm  \ref{tsm1.1-4}$\Rightarrow$\ref{tsm1.1-3}
$\Rightarrow$\ref{tsm1.1-1}} and  all four conditions are equivalent.
\end{thm}

Since $\mu$ is a positive  increasing function with $\mu(0_+)=0$ the finiteness restriction (\ref{eubc1.2}) is a condition on the growth $\mu$ at infinity, i.e.\ an implicit condition on  the possible growth of the coefficients of $H_\Omega$.
If the coefficients are uniformly bounded then $\mu(s)=O(s)$ as $s\to\infty$ and (\ref{eubc1.2}) is satisfied.
Then the  four conditions of the  theorem are equivalent.
 This retrieves the results of Theorems~1.2 and 1.3 of \cite{RSi4}.

The theorem is in part a restatement of standard results.
The equivalence of Conditions~\ref{tsm1.1-1}  and   \ref{tsm1.1-2} was established by  Davies
\cite{Dav14}, Theorem~2.2, whose arguments were based on earlier results of Azencott \cite{Aze}.
Although Davies assumptions were somewhat different his arguments apply with little modification
to the current setting.
The implication \ref{tsm1.1-2}$\Rightarrow$\ref{tsm1.1-3} is  a straightforward result which is established, for example, 
in \cite{Ebe} Lemma~1.16.
The implication   \ref{tsm1.1-3}$\Rightarrow$\ref{tsm1.1-4} follows as in the proof of Theorem~1.2 in  \cite{RSi4} for operators with $c_{ij}\in W^{1,\infty}(\Omega)$.

In the special case   $c(s)\sim s^2\,(\log s)^\alpha$ for large $s$ it follows that   $\mu(s)\sim (\log s)^{1-\alpha/2}$ and
 $\mu(s)\to\infty$ as $s\to\infty$ for  $\alpha\in[0,2\rangle$.
On the other hand if $\alpha\in [0,1]$ then (\ref{eubc1.2}) is satisfied for all sufficiently large $\lambda>0$.
Thus if $\alpha\in[0,2\rangle$ then  Markov uniqueness of $H_\Omega$ is equivalent to the capacity of the boundary being zero   and 
if $\alpha\in [0,1]$ then it is also equivalent to $L_1$-uniqueness of $H_\Omega$.

Note that if  $\Omega=\Ri^d$ the capacity condition is clearly satisfied and one concludes  that $H_{\Ri^d}$ is 
$L_1$-unique whenever (\ref{eubc1.2}) is satisfied for one large $\lambda>0$. 
 More generally we establish in Section~\ref{S4.2} that the   capacity condition  depends on the Hausdorff dimension of bounded subsets $A\subset \partial \Omega$  and  the order of degeneracy of $H_\Omega$ at~$A$.

\section{Submarkovian extensions}\label{S2}

The Friedrichs extension $H_{\Omega,D}$ of $H_\Omega$ is well known to be the largest submarkovian extension, i.e.\
the extension with the minimal form domain.
In this section we examine some basic properties of the smallest submarkovian extension, i.e.\ the extension with the maximal form domain (see, \cite{FOT} Section~3.3.3, \cite{Ebe} Section~3c or \cite{RSi4}, Section~3).
In particular we identify $H_{\Omega,N}$  as the smallest submarkovian extension.

We begin by discussing the imposition of  Neumann boundary conditions  on  a general submarkovian extension  $K_\Omega$  of $H_\Omega\in\ce_\Omega$.
Let $k_\Omega$ be the  Dirichlet form corresponding to $K_\Omega$.
Then    $D(k_\Omega)\cap L_\infty(\Omega)$ is an algebra.
Clearly $C_c^\infty(\Omega)\subseteq D(k_\Omega)\cap L_\infty(\Omega)$.
Thus one can define the truncated form $k_{\Omega,\chi}$ for each $\chi\in C_c^\infty(\Omega)$ by 
 $D(k_{\Omega,\chi})=D(k_\Omega)\cap L_\infty(\Omega)$
and 
\begin{equation}
k_{\Omega,\chi}(\varphi)=
k_\Omega(\varphi,\chi\varphi)-2^{-1}k_\Omega(\chi,\varphi^2)
\label{emin1.20}
\end{equation}
for $\varphi\in D(k_{\Omega,\chi})$.
The $k_{\Omega,\chi}$ have many properties similar to those of the forms $ h_{\Omega,\chi}$.
In particular the $k_{\Omega,\chi}$ are Markovian forms satisfying $0\leq k_{\Omega,\chi}\leq k_\Omega$.
Moreover,  if $\chi_1,\chi_2\in \cc_\Omega$ and $\chi_1\leq \chi_2$ then $k_{\Omega, \chi_1}\leq k_{\Omega, \chi_2}$
(see \cite{BH}, Proposition~I.4.1.1).
But it is not  evident that the $k_{\Omega,\chi}$ are closable.
This, however, is part of our first result.

\begin{thm}\label{text2.1}
Let  $H_\Omega\in \ce_\Omega$.
Further let $K_\Omega$ be a submarkovian extension of $H_\Omega$ and $k_\Omega$ the corresponding Dirichlet form.

If $\chi\in C_c^\infty(\Omega)$ then the truncated form $k_{\Omega,\chi}$ defined by $(\ref{emin1.20})$  is closable and the closure ${\overline k}_{\Omega,\chi}$ satisfies
${\overline k}_{\Omega,\chi}={\overline h}_{\Omega,\chi}$.
Therefore 
\[
h_{\Omega,N}\leq k_\Omega\leq h_{\Omega,D}\;.
\]
In particular $H_\Omega$ is Markov unique if and only if $h_{\Omega,N}=h_{\Omega,D}$.
\end{thm}
\proof\
The first step in the proof is a regularity property
 which extends a similar result for operators with bounded coefficients given by Theorem~1.1.IV in \cite{RSi4}.

\begin{lemma}\label{lmin1.1}
Let $K_\Omega$ be a  positive, self-adjoint extension of $H_\Omega$.
Then
\[
C_c^\infty(\Omega)D(K_\Omega)\subseteq D(\overline H_\Omega)\;.
\]
\end{lemma}
\proof\
First,  if $K_\Omega$ is a self-adjoint extension of $H_\Omega$ then $ H_\Omega\subseteq K_\Omega\subseteq H_\Omega^*$.
Therefore it suffices to establish  that $C_c^\infty(\Omega)D(H_\Omega^*)\subseteq D(\overline H_\Omega)$.
This property was proved for operators with bounded Lipschitz coefficients in Theorem~2.1 of \cite{RSi4}
but the proof is also valid for  operators with coefficients which are only locally bounded.
For example, if $\eta\in C_c^\infty(\Omega)$ with $\supp\eta=K$ and $V$ is a relatively compact subset of $\Omega$ with 
$K\subset V$ then to deduce that $\eta \,D(H_\Omega^*)\subseteq D(\overline H_\Omega)$ it suffices to prove that $\eta \,D(H_V^*)\subseteq D(\overline H_V)$ where $H_V$ is the restriction of $H_\Omega$ to $C_c^\infty(V)$.
Since the coefficients of $H_V$ are uniformly bounded  the result follows 
from Theorem~2.1 of \cite{RSi4}.~\hfill$\Box$

\bigskip

Next we prove the first statement of Theorem~\ref{text2.1}.

\begin{lemma} \label{lnsm2.2}
If $\chi\in C_c^\infty(\Omega)$ then $k_{\Omega,\chi}$ is closable and the closure ${\overline k}_{\Omega,\chi}$ satisfies
${\overline k}_{\Omega,\chi}={\overline h}_{\Omega,\chi}$.
\end{lemma}
\proof\
First $C_c^\infty(\Omega)D(K_\Omega)\subseteq D(\overline H_\Omega)$ by Lemma~\ref{lmin1.1}.
Now fix $\varphi\in D(K_\Omega)\cap L_\infty(\Omega)$.
Then for each $\chi\in \cc_\Omega$ one has $\chi\varphi\in D(\overline H_\Omega)$.
Moreover,
\[
k_\Omega(\varphi,\chi\varphi)=(K_\Omega\varphi, \chi\varphi)=(\varphi,\overline H_\Omega\chi\varphi)
\]
and 
\[
k_\Omega(\chi,\varphi^2)=(K_\Omega\chi,\varphi^2)=(H_\Omega\chi,\varphi^2)
\;.
\]
Therefore
\[
k_{\Omega,\chi}(\varphi)=(\varphi,\overline H_\Omega\chi\varphi)-2^{-1}(H_\Omega\chi,\varphi^2)
\;.
\]
Next choose a $\chi_1\in \cc_\Omega$ with $\chi_1=1$ on $\supp\chi$ and set 
$\varphi_1=\chi_1\varphi$.
It follows from Lemma~\ref{lmin1.1} that  $\varphi_1\in D(\overline H_\Omega)\cap L_\infty(\Omega)$.
Moreover,
\[
k_{\Omega,\chi}(\varphi)=(\varphi,\overline H_\Omega\chi\varphi_1)-2^{-1}(H_\Omega\chi,\varphi_1^2)
=(\varphi_1,\overline H_\Omega\chi\varphi_1)-2^{-1}(H_\Omega\chi,\varphi_1^2)=\overline h_{\Omega,\chi}(\varphi_1)
\;.
\]
The first equality  is obvious since $\supp\overline H_\Omega\chi=\supp\chi$.
The second equality  follows by approximating $\varphi$ 
in $L_2(\Omega)$ by a sequence $\varphi_n\in C_c^\infty(\Omega)$
and noting that 
\[
(\varphi_n,\overline H_\Omega\chi\varphi_1)=(H_\Omega\varphi_n,\chi\varphi_1)
=(H_\Omega\chi_1\varphi_n,\chi\varphi_1)=(\chi_1\varphi_n,\overline H_\Omega\chi\varphi_1)
\;.
\]
The third equality is also obvious.
But for $\chi$ and $\varphi$ fixed $\overline h_{\Omega,\chi}(\chi_1\varphi)$ is independent of the choice
of $\chi_1$.
Moreover, if $\chi_2$ is a second choice, with $\chi_2=1$ on $\supp\chi$ then
$\chi_1-\chi_2=0$ on $\supp \chi$ and $\overline h_{\Omega,\chi}((\chi_1-\chi_2)\varphi)=0$.
Therefore if $\chi_1\leq \chi_2\leq \ldots\leq \one_\Omega$ is an increasing family 
of $C_c^\infty$-functions with $\chi_n=1$ on $\supp\chi$ then $\overline h_{\Omega,\chi}((\chi_n-\chi_m)\varphi)=0$ but $\|\chi_n\varphi-\varphi\|_2\to0$.
This establishes that $\varphi\in D(\overline h_{\Omega,\chi})$ and $\overline h_{\Omega,\chi}(\varphi)=k_{\Omega,\chi}(\varphi)$.
Then, however,
\[
\overline  h_{\Omega,\chi}(\varphi)=k_{\Omega,\chi}(\varphi)\leq k_\Omega(\varphi)
\]
for all $\varphi\in D(K_\Omega)\cap L_\infty(\Omega)$.
Since $D(K_\Omega)$ is a core of $k_\Omega$ it follows by continuity that 
$\overline  h_{\Omega,\chi}(\varphi)=k_{\Omega,\chi}(\varphi)$ for all 
$\varphi\in D(k_\Omega)\cap L_\infty(\Omega)=D(k_{\Omega,\chi})$.
Therefore $\overline h_{\Omega,\chi}$ is a closed extension of $k_{\Omega,\chi}$.
Thus $k_{\Omega,\chi}$ is closable and its closure ${\overline k}_{\Omega,\chi}\subseteq {\overline h}_{\Omega,\chi}$.

But $h_{\Omega,\chi}(\psi)=k_{\Omega,\chi}(\psi)$ for all $\psi\in C_c^\infty(\Omega)$ and $C_c^\infty(\Omega)$
is a core of ${\overline h}_{\Omega,\chi}$ by definition.
Therefore ${\overline k}_{\Omega,\chi}\supseteq {\overline h}_{\Omega,\chi}$.
Combination of these conclusions gives ${\overline k}_{\Omega,\chi}={\overline h}_{\Omega,\chi}$.
\hfill$\Box$

\bigskip

One can now immediately deduce Theorem~\ref{text2.1}.
\smallskip

\noindent{\bf Proof of Theorem~\ref{text2.1}}
The first statement of the theorem  has been established by Lemma~\ref{lnsm2.2}.
Hence 
 \[
h_{\Omega,N}=\sup_{\chi\in \cc_\Omega} \overline k_{\Omega,\chi}\;.
 \]
 But $k_{\Omega,\chi}\leq k_\Omega$ for all $\chi\in\cc_\Omega$.
 Therefore $h_{\Omega,N}\leq k_\Omega$.

Finally $k_\Omega\supseteq h_\Omega$.
 Hence $k_\Omega\leq \overline h_\Omega=h_{\Omega,D}$.
 \hfill$\Box$
 
 \bigskip

The form $h_{\Omega,N}$ possesses a  {\it carr\'e du champ}  in the sense of  \cite{BH}, Section~I.4.
This is initially  defined  as  the bilinear form from $W^{1,2}_{\rm loc}(\Omega)\times W^{1,2}_{\rm loc}(\Omega)$ into $L_{1, \rm loc}(\Omega)$ given by

\[
\Gamma(\varphi\,;\psi)(x)=\sum^d_{i,j=1}c_{ij}(x)(\partial_i\varphi)(x)(\partial_j\psi)(x)
\]
and  $\Gamma(\varphi)=\Gamma(\varphi\,;\varphi)$.
Then 
\[
D(h_{\Omega,N})=\{\varphi\in W^{1,2}_{\rm loc}(\Omega): \sup_V\int_V dx\,\Gamma(\varphi)(x)<\infty\}
\]
where the supremum is over the relatively compact subsets $V $ of $ \Omega$ and
 \[
 h_{\Omega,N}(\varphi)=\sup_V\int_V dx\,\Gamma(\varphi)(x)
 \]
for all  $\varphi\in D(h_{\Omega,N})$.
It follows readily that if $\varphi\in D(h_{\Omega,N})$ then $\Gamma(\varphi)$ is a positive $L_1(\Omega)$-function
with $\|\Gamma(\varphi)\|_1=h_{\Omega,N}(\varphi)$.
The foregoing explicit identification of the form of the  minimal extension has been used in previous discussions
of Markov uniqueness,  \cite{FOT} Section~3.3.3, \cite{Ebe} Section~3c or \cite{RSi4}, Section~3.

A number of properties of  general submarkovian extension follows from the identification of the minimal extension.
If $k_\Omega$ is the form of the submarkovian extension $K_\Omega$ of $H_\Omega$ it follows from
Theorem~\ref{text2.1} that $D(k_\Omega)\subseteq D(h_{\Omega,N})$.
Therefore 
$k_\Omega$ possesses a  {\it carr\'e du champ} since $\Gamma(\varphi)\in L_1(\Omega)$ for all $\varphi\in D(k_\Omega)$. Moreover,  $k_\Omega(\varphi)=\|\Gamma(\varphi)\|_1$ for all $\varphi\in D(k_\Omega)$.
Further the form $h_{\Omega, N}$ is strongly local in the sense of \cite{FOT} and hence the restriction $k_\Omega$ is also strongly local.

\smallskip

Subsequently we need two  Dirichlet form implications of the elliptic regularity property.

\begin{cor}\label{cmin1.1}
Let $K_\Omega$ be a  submarkovian  extension of $H_\Omega\in \ce_\Omega$ and $k_\Omega$ the corresponding
Dirichlet form.
Then
\[
C_c^\infty(\Omega)D(k_\Omega)\subseteq D(\overline h_\Omega)\;.
\]
\end{cor}
\proof\
If $\eta\in C_c^\infty(\Omega)$ and $\varphi\in D(K_\Omega)$ then it follows from Lemma~\ref{lmin1.1} that 
$\eta\,\varphi\in D(\overline H_\Omega)\subset D(K_\Omega)\subseteq D(k_\Omega)\subseteq D(h_{\Omega,N})$.
Moreover,
\begin{eqnarray*}
\overline h_\Omega(\eta\,\varphi)=h_{\Omega,N}(\eta\,\varphi)\leq 
2\int_\Omega\Gamma(\eta)\,\varphi^2+2\int_\Omega\eta^2\,\Gamma(\varphi)
\leq 2\,\left(\|\Gamma(\eta)\|_\infty+\|\eta\|_\infty^2\right)\|\varphi\|_{D(k_\Omega)}^2
\;.
\end{eqnarray*}
Since $D(K_\Omega)$ is a core of $k_\Omega$ with respect to the $D(k_\Omega)$-graph norm this estimate extends to all $\varphi\in D(k_\Omega)$ by continuity.
The 
statement of the corollary follows immediately.
\hfill$\Box$

\bigskip

\begin{cor}\label{cmin1.2}
If $H_\Omega\in \ce_\Omega$ then
$C_c^\infty(\Ri^d)D( h_{\Omega,N})\subseteq D( h_{\Omega,N})$.
\end{cor}
\proof\
Fix  $\rho\in C_c^\infty(\Ri^d)$ and $\chi\in\cc_\Omega$.
If  $\varphi\in C_c^\infty(\Omega)$.
Then $\rho\,\varphi\in C_c^\infty(\Omega)$ and $\|\rho\,\varphi\|_2\leq \|\rho\|_\infty\|\varphi\|_2$.
Moreover,
\[
h_{\Omega, \chi}(\rho\,\varphi)\leq 2\,\|\rho\|_\infty^2\,h_{\Omega, \chi}(\varphi)+2\,a_\rho\,\|\nabla\rho\|_\infty^2\,\|\varphi\|_2^2
\]
where $a_\rho=\sup_{x\in\supp\rho}\|C(x)\|$.
Therefore, by continuity, $\rho\, D(\overline h_{\Omega,\chi})\subseteq D(\overline h_{\Omega,\chi})$ and 
\[
\|\rho\,\varphi\|_{D(\overline h_{\Omega,\chi})}\leq a(\rho)\,\|\varphi\|_{D(\overline h_{\Omega,\chi})}
\]
for all $\varphi\in D(\overline h_{\Omega, \chi})$ with $a(\rho)=2\,(a_\rho\,\|\nabla\rho\|_\infty^2+\|\rho\|_\infty^2)$.
Since this estimate is uniform for $\chi\in \cc_\Omega$ it follows that $\rho\, D(h_{\Omega, N})\subseteq D(h_{\Omega, N})$
and $\|\rho\,\varphi\|_{D(h_{\Omega, N})}\leq a(\rho)\,\|\varphi\|_{D(h_{\Omega, N})}$.
\hfill$\Box$

\section{$L_1$-uniqueness}\label{S3}

In this section we prove Theorem~\ref{tsm1.1}.
Much of the proof consists of refinements of previous arguments.
\smallskip

\noindent \ref{tsm1.1-1}$\Leftrightarrow$\ref{tsm1.1-2}$\;$
This equivalence was established by  Davies
\cite{Dav14}, Theorem~2.2, for a large class of second-order elliptic operators with
smooth coefficients.
But his arguments extend to the current situation with only minor modifications.
We omit further details

\smallskip

\noindent \ref{tsm1.1-2}$\Rightarrow$\ref{tsm1.1-3}$\;$
This is a general structural result which is proved, for example, in Lemma~1.1.6 of \cite{Ebe}.

\smallskip

\noindent\ref{tsm1.1-3}$\Rightarrow$\ref{tsm1.1-4}$\;$
First note that Markov uniqueness of $H_\Omega$ is equivalent to the identity $h_{\Omega,N}=h_{\Omega,D}$ by Theorem~\ref{text2.1}.
But in general $h_{\Omega,N}\supseteq h_{\Omega,D}$
and $C_c^\infty(\Omega)$ is a core of $h_{\Omega,D}$.
Therefore Markov uniqueness of $H_\Omega$ implies  that 
$C_c^\infty(\Omega)$ is  a core of $h_{\Omega,N}$.

Secondly, let $\psi\in D(h_{\Omega,N})\cap L_\infty(\Omega)$ with $\psi=1$ on $U\cap\Omega$ where
$U$ is an open subset containing $\partial\Omega$.
Then since $C_c^\infty(\Omega)$ is a core of $h_{\Omega,N}$  there is a sequence $\psi_n\in C_c^\infty(\Omega)$ such that $\|\psi-\psi_n\|_{D(h_{\Omega,N})}\to0$ as $n\to\infty$.
Set $\varphi_n=\psi-\psi_n$.
Then $\varphi_n\in D(h_{\Omega,N})$, $\|\varphi_n\|_{D(h_{\Omega,N})}\to0$ and 
since $\psi_n$ has compact support there is an open subset $U_n$ containing $\partial\Omega$ such that $\varphi_n=1$
on $(U\cap U_n)\cap \Omega$.
Therefore $\capp_\Omega(\partial\Omega)=0$.

\medskip

Combination of the foregoing observations  establishes the first statement of Theorem~\ref{tsm1.1}.
Now we turn to the proof of the second statement.
\medskip

\noindent \ref{tsm1.1-4}$\Rightarrow$\ref{tsm1.1-3}$\;$
Assume
$\mu(s)\to\infty$ as $s\to\infty$ where $\mu $ is  defined by (\ref{eubc1.10}) and (\ref{eubc1.11}).
Then  $H_\Omega$ is Markov unique if and only if 
$C_c^\infty(\Omega)$ is a core of $h_{\Omega,N}$.
Thus it is necessary to   demonstrate that 
each $\varphi\in D(h_{\Omega,N})\cap L_\infty(\Omega)$ can be approximated
in the $D(h_{\Omega,N})$-graph norm by a sequence $\varphi_n\in C_c^\infty(\Omega)$.

Fix $\varphi\in D(h_{\Omega,N})\cap L_\infty(\Omega)$.
Next fix  $\rho\in C^\infty_c(\Ri)$ with $0\leq\rho\leq 1$, $\rho(s)=1$ if $s\leq 1$ and $\rho(s)=0$ if $s\geq 2$.
Then define $\rho_n$ by $\rho_n(x)=\rho(n^{-1}\mu(|x|))$.
It follows that  $\rho_n(x)=1$ if $\mu(|x|)\leq n$ and $\rho(x)=0$ if $\mu(|x|)\geq 2n$.
Moreover,  $\|\Gamma(\rho_n)\|_\infty\leq b^2\,n^{-2}$
with $b=\|\rho'\|_\infty$.
Then $\rho_n\,\varphi\in D(h_{\Omega,N})\cap L_\infty(\Omega)$ by Corollary~\ref{cmin1.2} and
\begin{eqnarray*}
\|\varphi-\rho_n\,\varphi\|_{D(h_{\Omega,N})}^2&\leq& 2\int_\Omega\varphi^2\,\Gamma(\rho_n)
+2\int_\Omega(\one_\Omega-\rho_n)^2\,\Gamma(\varphi)+\|(\one_\Omega-\rho_n)\varphi\|_2^2\\[5pt]
&\leq &2\,b^2\,n^{-2}\|\varphi\|_2^2+\int_\Omega(\one_\Omega-\rho_n)^2(2\,\Gamma(\varphi)+\varphi^2)
\;.
\end{eqnarray*}
The first term on the right hand side clearly  tends to zero as $n\to\infty$.
But  it follows by construction that $(\one_\Omega-\rho_n)^2\to0$ pointwise on $\Omega$. 
Therefore  the second term also tends to zero by the Lebesgue dominated convergence theorem.
Thus $\varphi$ is approximated by the sequence $\rho_n\varphi$ in the $D(h_{\Omega,N})$-graph norm.

Next since $\capp_\Omega(\partial\Omega)=0$ one may choose $\chi_n\in D(h_{\Omega,N})$ and
open subsets $U_n\supset\partial\Omega$ such that $0\leq\chi_n\leq1$, $\|\chi_n\|_{D(h_{\Omega,N})}\leq n^{-1}$
and $\chi_n\geq1$ on $U_n\cap\Omega$.
But since $h_{\Omega,N}$ is a Dirichlet form one may assume $\chi_n=1$ on  $U_n\cap\Omega$.
Then with  $\varphi_n=(\one_\Omega-\chi_n)\rho_n\varphi$  one has
\[
\lim_{n\to\infty}\|\varphi-\varphi_n\|_{D(h_{\Omega,N})}\leq \lim_{n\to\infty}\|\chi_n\rho_n\varphi\|_{D(h_{\Omega,N})}
\]
by the Cauchy--Schwarz estimate and the conclusion of the previous paragraph.
But
\[
\|\chi_n\rho_n\varphi\|_{D(h_{\Omega,N})}^2=h_{\Omega,N}(\chi_n\rho_n\varphi)+\|\chi_n\rho_n\varphi\|_2^2
\]
and the second term on the right hand side tends to zero as $n\to\infty$ because $\|\chi_n\rho_n\varphi\|_2\leq \|\chi_n\|_2\|\varphi\|_\infty$.
The first term on the right can, however, be estimated by
\begin{eqnarray*}
h_{\Omega,N}(\chi_n\rho_n\varphi)&\leq&2\int_\Omega\varphi^2\,\Gamma(\chi_n)
+4\int_\Omega\varphi^2\,\Gamma(\rho_n)+4\int_\Omega\chi_n^2\,\Gamma(\varphi)\\[5pt]
&\leq&2\,\|\varphi\|_\infty^2\,h_{\Omega,N}(\chi_n)+4\,b^2\,n^{-2}\,\|\varphi\|_2^2
+4\int_\Omega\chi_n^2\,\Gamma(\varphi)
\end{eqnarray*}
since $\|\Gamma(\rho_n)\|_\infty\leq b^2\,n^{-2}$.
The first term on the right hand side tends to zero because $h_{\Omega,N}(\chi_n)\leq n^{-1}$
and the second obviously  tends to zero.
Finally  the third tends to zero by an equicontinuity estimate because $\chi_n^2\leq 1$ and  $\Gamma(\varphi)\in L_1(\Omega)$.
Thus one now concludes that $\varphi$ is approximated by the sequence $\varphi_n$ in the $D(h_{\Omega,N})$-graph norm.

Finally, $\supp\varphi_n\subseteq\Omega_n=((\supp\rho_n)\cap\Omega)\cap(\Omega\backslash(U_n\cap\Omega))$
and $\Omega_n$ is a relatively compact subset of $\Omega$. 
Therefore $H_{\Omega, N}$ is strongly elliptic in restriction to $\Omega_n$.
Consequently
$\varphi_n$, and hence $\varphi$, can be approximated by a sequence of $C_c^\infty(\Omega_n)$-functions in the $D(h_{\Omega,N})$-graph norm.

\smallskip

This completes the proof of the second statement of Theorem~\ref{tsm1.1}.
Now we turn to the proof of the third statement.
By the foregoing  it suffices to prove the following.

\smallskip

\noindent \ref{tsm1.1-3}$\Rightarrow$\ref{tsm1.1-1}$\;$
The proof is an elaboration  of  the argument used to demonstrate the comparable implication in  Theorem~1.3 in \cite{RSi4}.

\begin{prop}\label{pubc3.1}
Assume $H_\Omega\in\ce_\Omega$ is Markov unique.
Further assume that 
\[
\int_0^\infty ds\, s^{d/2}\, e^{-\lambda\,\mu(s)^2}<\infty
\]
for one  $\lambda>0$
where 
$\mu(s)=\int^s_0c^{-1/2}$  with $c$ defined by $(\ref{eubc1.10})$.

\smallskip

Then $H_\Omega$ is conservative.
\end{prop}
\proof\
The proof is in several steps.

\smallskip

\noindent{\bf Step~1: $\;\Omega$ bounded.}$\;$
If $\Omega$ is bounded then $H_\Omega$ is conservative by  Step~1 in the proof of Theorem~1.3 in \cite{RSi4}.
 Therefore we now assume that $\Omega$ is unbounded.

\smallskip

\noindent{\bf Step~2: Bounded approximation.}$\;$
The second step consists of introducing an increasing sequence of bounded sets $\Omega_n$ and conservative
operators   $H_{\Omega;n}\in \ce_{\Omega_n}$ which approximate $H_\Omega$ in a suitable manner.

First, fix  $\rho\in C^\infty_c(\Ri)$ with $0\leq\rho\leq 1$, $\rho(s)=1$ if $|s|\leq 1$ and $\rho(s)=0$ if $|s|\geq 2$.
Then  introduce the sequence $\rho_n$ by $\rho_n(x)=\rho(n^{-1}|x|)$.
Thus $\rho_n(x)=1$ if $|x|\leq n$ and $\rho(x)=0$ if $|x|\geq 2n$.
The family of functions $\rho_n$ is monotonically increasing.
Set $B_n=\{x\in\Ri^d:|x|<n\}$ and $\Omega_n=\Omega\cap B_{2n}$.
Note that $\Omega_n$ is bounded.

Secondly, define $H_{\Omega; n}\in \ce_{\Omega_n}$ as the operator with coefficients $\rho_n\,c_{ij}$ acting on $L_2(\Omega_n)$.
Then it follows that $H_{\Omega;n}$ is Markov unique since the capacity of $\partial \Omega_n$ with respect to the 
Neumann  form associated with  $H_{\Omega; n}$ is zero.
Therefore $H_{\Omega; n}$ is conservative by Step~1.
Then if $H_n$ is the  extension to $L_2(\Omega)$ of the 
 unique submarkovian extension $H_{\Omega; n, N}(=H_{\Omega; n, D})$ of $H_{\Omega; n}$ acting on
 $L_2(\Omega_n)$, i.e.\ if
 $H_n=H_{\Omega;n,N}\oplus 0$ with $L_2(\Omega)=L_2(\Omega_n)\oplus L_2(\Omega_n)^\perp$,
then $H_n$ is conservative.

\smallskip

\noindent{\bf Step~3: $L_2$-convergence.}$\;$
The third step is to establish strong convergence  on $L_2(\Omega)$ of the  semigroups $S^{(n)}$ generated by the $H_n$
to the semigroup $S$ generated by the unique submarkovian extension $H_{\Omega,N}(=H_{\Omega, D})$ of $H_\Omega$.
This follows by a monotone convergence argument.
The closed form $h_n$ corresponding to $H_n$ on $L_2(\Omega)$ is given by 
$ h_n(\varphi)=h_{\Omega;n, N}(\varphi)$ for all $\varphi\in C_c^\infty(\Omega)$ and then by closure for all $\varphi\in D(h_n)$.
Since the $\rho_n$ are a monotonically increasing family of functions on $\Ri^d$ the  forms $ h_n$  are a monotonically increasing family of Dirichlet forms. 
If $h=\sup_{n\geq1}h_n$ then $h$ is a Dirichlet form.

It follows from the monotonic increase of the forms $h_n$ that the operators $H_n$ converge in the strong resolvent sense on $L_2(\Omega)$ to the operator  $H$ corresponding to $h$ (see, for example, \cite{Kat1}, Section~VIII.4, or \cite{MR}, Section~I.3).
Moreover,  the semigroups $S^{(n)}$ converge strongly  on $L_2(\Omega)$ to the submarkovian semigroup $S$ generated by  $H$.
It also  follows readily that  $H$ is a submarkovian extension of $H_\Omega$.
Therefore  $H=H_{\Omega, N}(=H_{\Omega,D})$,  by Markov uniqueness.

\medskip

Our next aim is to prove that the semigroups $S^{(n)}$ converge strongly   to  $S$ on $L_1(\Omega)$.
Following a tactic used in \cite{RSi2} \cite{RSi4}, we convert the  $L_2$-convergence of the semigroups into $L_1$-convergence by the use of suitable off-diagonal bounds.

\medskip

\noindent{\bf Step~4: $L_2$-off-diagonal bounds.}$\;$
Let
\begin{equation}
D_n=\{\psi\in W^{1,\infty}(\Omega):\sum^d_{i,j=1}\rho_n\,c_{ij}\,(\partial_i\psi)(\partial_j\psi)\leq1\}
\;.
\label{ecsg4.0}
\end{equation}
The corresponding   Riemannian (pseudo-)distance is defined by
\begin{equation}
d_n(x\,;y)=\sup_{\psi\in D_n}\,(\psi(x)-\psi(y))
\label{ecsg4.1}
\end{equation}
for all $x,y\in\Omega$.
This function has the metric properties of a distance but it takes the value infinity if either $x$ or $y$ is not in $\Omega_n$.
Secondly, introduce the corresponding set-theoretic distance by
\[
d_n(A\,;B)=\inf_{x\in A,\,y\in B}d_n(x\,;y)
\]
where $A$ and $B$ are general measurable subsets of $\Omega$.
Finally define  $D$  by setting $\rho_n=\one_\Omega$ in (\ref{ecsg4.0}). 
Then $D\subseteq D_n$ and the corresponding Riemannian distance  $d(\,\cdot\,;\,\cdot\,)$, defined  in analogy with (\ref {ecsg4.1}), satisfies $d(x\,;y)\leq d_n(x\,;y)$.

\begin{lemma}\label{lclo4}
If $A,B$  are  open subsets of $\Omega$ then 
\[
\sup_{n\geq1}|(\varphi_A, S^{(n)}_t\varphi_B)|\vee |(\varphi_A, S_t\varphi_B)|
\leq e^{-d(A;B)^2(4t)^{-1}}\|\varphi_A\|_2\,\|\varphi_B\|_2
   s\]
for all $\varphi_A\in L_2(A)$, $\varphi_B\in L_2(B)$ and $t>0$ with the convention
$e^{-\infty}=0$.
\end{lemma}

Bounds of this type have now been derived by many authors 
(see, for example, \cite{Aus2} \cite{CGT}  \cite{Dav12} \cite{Gri3}  \cite{Stu4} \cite{Stu2}) under a  variety of  ellipticity assumptions.
A proof  applicable in the current context can be found in \cite{RSi2}, Section~4.
The bounds for $S^{(n)}$ are initially in terms of $d_n(A\,;B)$ but $d_n(A\,;B)\leq d(A\,;B)$.
Then since the $S^{(n)}_t$ are $L_2$-convergent to $S_t$ the bounds also hold for $S$.
\smallskip

Next  $C(x)\leq c(|x|)\,I$ for all $x\in\Omega$.
Therefore
 \begin{equation}
D_n\supseteq\widehat D_n= \{\psi\in W^{1,\infty}(\Omega):\rho(n^{-1}|x|)\,c(|x|)\,|(\nabla\psi)(x)|^2\leq1\}
\;.
\label{eubc3.21}
\end{equation}
Consequently
 \begin{equation}
\hat d_n(x\,;y)=\sup_{\psi\in \widehat D_n}\,(\psi(x)-\psi(y))\leq d_n(x\,;y)
\label{eubc3.22}
\end{equation}
for all $x,y\in\Omega$.
Moreover, if   $\widehat D$  is defined  by setting $\rho=\one_\Omega$ in (\ref{eubc3.21}) and $\hat d(\,\cdot\,;\,\cdot\,)$
is  defined  in analogue  with (\ref {eubc3.22}) then  $\hat d(x\,;y)\leq d(x\,;y)$ for all $x,y\in\Omega$.
Thus the bounds of Lemma~\ref{lclo4} are also valid with $d(A\,;B)$ replaced by $\hat d(A\,;B)$.

If $A, B$ are bounded open sets with $ A\subset \Omega\cap B_m$ and $B\subset \Omega\cap (B_M)^{\rm c}$ where $M>m\geq1$ define
\[
\rho_m=\sup_{x\in \Omega\cap B_m}\hat  d(x\,;0)\;\;\;\;\;{\rm and}\;\;\;\;\;\nu_M=\inf_{x\in \Omega\cap (B_M)^{\rm c}}\hat d(x\,;0)
\,.
\]
Then it follows from the triangle inequality $\hat d(x\,;0)\leq \hat d(x\,;y)+\hat d(y\,;0)$ that
\[
\nu_M\leq \inf_{x\in  B}\hat d(x\,;0)
\leq\inf_{x\in  B}\hat d(x\,;y)+\hat d(y\,;0)
\leq \inf_{x\in  B}\hat d(x\,;y)+\rho_m
\]
for all $y\in A$.
Therefore
\[
\hat d(A\,;B)\geq \nu_M-\rho_m\geq0
\]
where the last inequality follows because $M>m$.
But it follows directly from the definition of $\hat d(\,\cdot\,;\,\cdot\,)$ that
\[
\hat d(x\,;0)=\int^{|x|}_0ds\,c(s)^{-1/2}
\]
for all $x\in\Omega$.
Therefore
\[
\rho_m=\int^m_0ds\,c(s)^{-1/2}=\mu(m)\;\;\;\;\;{\rm and}\;\;\;\;\;\hat \nu_M=\int^M_0ds\,c(s)^{-1/2}=\mu(M)
\,.
\]
Hence \[
\hat d(A\,;B)\geq\mu(M)-\mu(m)\geq0
\;.
\]
Consequently one has the following variation of Lemma~\ref{lclo4}.

\begin{lemma}\label{lclo4.1}
If $M>m\geq1$ and $A,B$ are  bounded  open sets  with $ A\subset \Omega\cap B_m$ and $B\subset \Omega\cap (B_M)^{\rm c}$ then 
\[
\sup_{n\geq 1}|(\varphi_A, S^{(n)}_t\varphi_B)|\vee |(\varphi_A, S_t\varphi_B)|
\leq e^{\mu(m)^2(4t)^{-1}} e^{-\mu(M)^2(8t)^{-1}}\|\varphi_A\|_2\,\|\varphi_B\|_2
\]
for all  $\varphi_A\in L_2(A)$, $\varphi_B\in L_2(B)$ and $t>0$.
\end{lemma}
\proof\
The bounds on $|(\varphi_A, S_t\varphi_B)|$ follow directly from the bounds of Lemma~\ref{lclo4},
the foregoing observation that $\hat d(A\,;B)\geq\mu(M)-\mu(m)\geq0$ and the estimate
\[
 (\mu(M)-\mu(m))^2\geq 2^{-1}\mu(M)^2-\mu(m)^2
 \;.
 \]
The bounds on $\sup_{n\geq 1}|(\varphi_A, S^{(n)}_t\varphi_B)|$ follow by similar reasoning since $\hat d_n(A\,;B)\geq \hat d(A\,;B)$.
We omit further details.\hfill$\Box$

\bigskip

Now we are prepared for the key estimate.

\begin{lemma}\label{lubc5}
There is a $b>0$ such that if  $M>m\geq1$ then
\[
\sup_{n\geq1}|(\one_{(B_M)^{\rm c}}, S^{(n)}_t\varphi)|\vee
|(\one_{(B_M)^{\rm c}}, S_t\varphi)|\leq b\,e^{\mu(m)^2(4t)^{-1}}\int^\infty_M ds\,s^{d/2}\,e^{-\mu(s)^2(8t)^{-1}}\,\|\varphi\|_2
\]
for all  $\varphi\in L_2(\Omega\cap B_m)$ and  $t>0$.
\end{lemma}
\proof\
The proof is a variation of   an argument of \cite{ERSZ1}.
Let $C_p=B_{p+1}\backslash B_p$.
It follows that  $(B_M)^{\rm c}=\bigcup_{p\geq M}C_p$.
If $A$ is a bounded open set with  $\supp\varphi\subset A\subseteq \Omega\cap B_m$ then
by  Lemma~\ref{lclo4.1} 
\[
|(\one_{(B_M)^{\rm c}}, S_t\varphi)|=|\sum_{p\geq M}(\one_{C_p}, S_t\varphi)|
\leq e^{\mu(m)^2(4t)^{-1}}\sum_{p\geq M}e^{-\mu(p)^2(8t)^{-1}}\,|B_{p+1}|^{1/2}\,\|\varphi\|_2
\;.
\]
But the sum is a Riemann approximation to the integral occurring in the statement of the lemma.
Therefore  the bounds for $|(\one_{(B_M)^{\rm c}}, S_t\varphi)|$ follow immediately.
The bounds for $|(\one_{(B_M)^{\rm c}}, S^{(n)}_t\varphi)|$ follow by similar reasoning.
\hfill$\Box$

\smallskip

\noindent{\bf Step~5:  $L_1$-convergence.}$\;$
The fifth step  consists of  proving that the semigroups $S^{(n)}_t$ are strongly convergent on $L_1(\Omega)$ to 
$S_t$ (see \cite{RSi2}, Proposition~6.2, for a similar result).

 Since the semigroups $S^{(n)}_t$ and $S_t$  are all submarkovian  it suffices to prove convergence
 on a  subset of $L_1(\Omega)$ whose span is dense.
 In particular it suffices to prove convergence on positive functions in $L_1(A)\cap L_2(A)$ for each bounded open
 subset $A$ of $\Omega$.
 
 Fix $A\subset\Omega\cap B_m$ and $\varphi\in L_1(A)\cap L_2(A)$.
 Assume $\varphi$ is positive.
Then
\begin{eqnarray*}
\|(S^{(n)}_t-S_t)\varphi\|_1&\leq &
\|\one_{B_M}(S^{(n)}_t-S_t)\varphi\|_1
+\|\one_{(B_M)^{\rm c}}S^{(n)}_t\varphi\|_1+\|\one_{(B_M)^{\rm c}}S_t\varphi\|_1\\[5pt]
&\leq&|B_M|^{1/2}\|(S^{(n)}_t-S_t)\varphi\|_2
+(\one_{(B_M)^{\rm c}},S^{(n)}_t\varphi)+(\one_{ (B_M)^{\rm c}},S_t\varphi)\\[5pt]
&\leq&|B_M|^{1/2}\|(S^{(n)}_t-S_t)\varphi\|_2+
2\,b\,e^{\mu(m)^2(4t)^{-1}}\int^\infty_M ds\,s^{d/2}\,e^{-\mu(s)^2(8t)^{-1}}\,\|\varphi\|_2
\end{eqnarray*}
where we have used the positivity of the semigroups and the functions to express the $L_1$-norms as pairings between $L_1$ and $L_\infty$.
The last step uses Lemma~\ref{lubc5}.
But the integral is convergent for one $t=t_0>0$, by assumption.
Therefore it is convergent for all  $t\in\langle0,t_0]$.
Then  since $S^{(n)}_t$ is $L_2$-convergent to $S_t$ for all $t>0$  and since the last term on the right hand side
 converges to zero as $M\to\infty$ for each $t\in\langle0,t_0]$ it follows that $S^{(n)}_t$ is $L_1$-convergent to $S_t$ for all $t\in\langle0,t_0]$.
Finally it follows from the semigroup property and contractivity that  $S^{(n)}_t$ is $L_1$-convergent to $S_t$ for all $t>0$.

\smallskip

\noindent{\bf Step~6:  Conservation.}$\;$
The conservation  property for $S$ now follows because the approximating semigroups $S^{(n)}$ are conservative, by Step~2, and are $L_1$-convergent to $S$, by Step~5.
Therefore
\[
(\one_\Omega, S_t\varphi)=\lim_{n\to\infty}(\one_\Omega,S^{(n)}_t\varphi)
=\lim_{n\to\infty}(S^{(n)}_t\one_\Omega,\varphi)=(\one_\Omega, \varphi)
\]
for all $\varphi\in L_1(\Omega)$.
Hence $S_t\one_\Omega=\one_\Omega$.  
\hfill$\Box$

\bigskip

This completes the proof of Proposition~\ref{pubc3.1}  and the third statement of Theorem~\ref{tsm1.1}.

\section{Illustrations and examples}\label{S4}

In this section we illustrate the foregoing results with some applications and examples.

Theorem~\ref{tsm1.1} established that $L_1$-uniqueness of $H_\Omega\in \ce_\Omega$ is a consequence of two distinct properties,
a capacity condition on the boundary and a growth condition on the coefficients.
Therefore we separate the initial discussion into two parts each concentrating one of these conditions.

\vspace{-4mm}

\subsection{Growth properties}\label{S4.1}
If $\Omega=\Ri^d$ then the capacity condition plays no role and so we begin by considering this  case.
We continue to use the function $c$ and the corresponding function $\mu$ defined by (\ref{eubc1.10}) and 
(\ref{eubc1.11}), respectively, as measures of the coefficient growth.
The following statement combines  the $L_1$-properties which follow from the foregoing with the comparable 
$L_2$-properties established earlier by Davies et al. (see \cite{Dav14}  and references therein).

\begin{prop}\label{pubc4.1}
Let $H\in\ce_{\Ri^d}$.
Then the following are valid:
\begin{tabel}
\item\label{pubc4.1-1} 
If $\mu(s)\to\infty$ as $s\to\infty$ then
$H$ is $L_2$-unique,

\medskip

 \item\label{pubc4.1-2}
  If $\int_0^\infty ds\, s^{d/2}\, e^{-\lambda\,\mu(s)^2}<\infty $ for one  $\lambda>0$ then
$H$ is $L_1$-unique.
\end{tabel}
\end{prop}

The second statement is  a direct consequence of the second statement of Theorem~\ref{tsm1.1} since one automatically
has   $\capp_\Omega(\partial\Omega)=0$.
The first statement follows from \cite{Dav14}, Theorem~3.2.
 This theorem asserts that $H$ is essentially self-adjoint, i.e.\ $L_2$-unique, 
 if there exists a strictly positive differentiable function $\eta$ over $\Ri^d$
with $\|\Gamma(\eta)\|_\infty<\infty$ such that $\eta(x)\to\infty$ as $x\to\infty$.
But it follows by assumption that $\eta(x)=\int^{|x|}_0dt\,c(t)^{-1/2}$ satisfies these properties.

\medskip

Next we consider a special case of the growth property examined earlier by Davies \cite{Dav14} and Eberle \cite{Ebe}.

\begin{exam}\label{exubc4.1}
Assume that    
$c(s)\leq a\,s^2\,(\log s)^\alpha$ for some $a>0$, $\alpha\geq0$  and all large $s$. 
In this case  $\mu(s)\geq b\,(\log s)^{1-\alpha/2}$ with $b>0$ for all large~$s$.
Therefore  $\mu(s)\to\infty$ as $s\to\infty$ if $\alpha<2$ and $H$ is $L_2$-unique by the first statement of Proposition~\ref{pubc4.1}.
But Davies has demonstrated by specific example,   \cite{Dav14} Example~3.5, that if $d\geq2 $ and $\alpha>2$ then
$L_2$-uniqueness can fail (see also \cite{Ebe} Chapter~2, Section~c).
Finally Theorem~2.3 in \cite{Ebe} treats the borderline case $\alpha=2$. 
This theorem establishes that  if  $c(s)\leq a\,s^2\,(\log s)^2$ for all large $s$ then $H$ is not only $L_2$-unique but also $L_p$-unique for all $p\in\langle1,2]$.

Next if $\alpha\leq 1$ then  the  second statement of Proposition~\ref{pubc4.1} establishes that the  bound $H$ is  $L_1$-unique.
Indeed if $\alpha\leq 1$ then
  $\mu(s)\geq b\,(\log s)^{1-\alpha/2}$ for $s$ large
and the integral (\ref{eubc1.2})  is finite for  large $\lambda$.
Therefore Proposition~\ref{pubc4.1} establishes that $H$ is $L_1$-unique.
But $L_1$-uniqueness can fail if  
$c(s)\sim s^2\,(\log s)^\alpha $ with $\alpha>1$ for  large $s$.
To verify this let $d=1$ and define the positive  $L_\infty$-function $\psi$ on $\Ri $ by 
  \begin{equation}
  \psi(x)=1-a\,(\log (\log2))\,(\log(\log(2+|x|))^{-1}
 \label{eubc4.4}
 \end{equation}
   with $a\in\langle0,1\rangle$.
    Thus $\psi(0)=1-a>0$ and $\psi(x)\to1$ as $|x|\to\infty$. 
  Then define $c$  by $c(x)=\psi'(x)^{-1}\int^x_0ds\,\psi(s)$.
  It is evident that $c$ is strictly positive and $c\in W^{1,\infty}_{\rm loc}(\Ri)$. 
Moreover, $c(x)\sim  |x|^2(\log|x|)(\log(\log|x|))$ as $|x|\to\infty$.
But if $H$ is the corresponding operator on $C_c^\infty(\Ri)$, i.e.\ if $H\varphi= -(c\,\varphi')'$,
then $(I+H)\psi=\psi-(c\,\psi')'=0$.
Therefore the range of $I+H$ is not $L_1$-dense and 
  $H$ is not $L_1$-unique.
  
  Therefore within this class of examples the growth bound $c(s)\leq a\,s^2\,(\log s)$ is optimal for $L_1$-uniqueness
and the bound $c(s)\leq a\,s^2\,(\log s)^2$ is optimal for $L_2$-uniqueness.
\hfill$\Box$
\end{exam}
\smallskip

Note that if $d=1$,  $\Omega=\langle0,\infty\rangle$ and one repeats the foregoing construction with $\psi$ given by (\ref{eubc4.4})  but with $a=1$ then $c$ is strictly positive and $c(x)= O(x)$ as $x\to0$. 
  Moreover, $c(x)\sim x^2(\log x)(\log(\log x))$ for all large $x$.
Therefore the corresponding operator  $H$ is Markov unique, by \cite{RSi3} Theorem~2.7, but again it is not $L_1$-unique.
In fact it is not $L_p$-unique for any $p\in[1,\infty\rangle$.

\bigskip

The function $\mu$ is a lower bound on the Riemmanian distance to infinity measured with respect to 
the metric  $ C^{-1}$ associated with the operator $H_\Omega$.
  If one has more detailed information on the geometry one can obtain stronger conclusion by the same general
  reasoning. 
  This is illustrated by the following example of a Gru\v{s}in-type operator.

\smallskip

Let $d=2$ and $\Omega=\Omega_+\cup \Omega_-$ with $\Omega_\pm=\{x=(x_1,x_2): \pm x_1>0\}$.
Define the Gru\v{s}in operator $H$ by $D(H)=C_c^\infty(\Omega)$ and
\[
(H\varphi)(x)=-\partial_1(c_1(x_1)\partial_1\varphi)(x)-c_2(x_1)(\partial_2^2\varphi)(x)
\]
where $c_1, c_2\in W^{1,\infty}_{\rm loc}(\Ri\backslash\{0\})$ are strictly positive  and $c_i(x)\sim|x|^{(2\delta_i,2\delta'_i)}$
with $\delta_i,\delta'_i\geq0$.
Here we use the notation of \cite{RSi2a} \cite{RSi2}.
Specifically $s^{(\alpha, \alpha')}=s^\alpha$ if $\alpha\leq 1$ and  $s^{(\alpha, \alpha')}=s^{\alpha'}$ if $\alpha\geq 1$ and functions $f,g$ satisfy the relation $f\sim g$ if there are $a,a'>0$ such that $a\,f\leq g\leq a'\,f$.
We assume that $\delta_1,\delta_1'<1$ but there are no upper bounds on $\delta_2$ and $\delta_2'$.
Thus $H\in \ce_\Omega$ with
\[
C(x)=\left( \begin{array}{cc}
                     c_1(x_1)&0\\[5pt]
                     0&c_2(x_1) 
                     \end{array}\right)\;.
\]
Therefore $\|C(x)\|=c_1(x_1)\vee c_2(x_1)\leq a\,|x_1|^{2(\delta_1'\vee \delta_2')}\leq a\,|x_1|^{2(1\vee \delta_2')}$ for all $|x_1|\geq1$.
Although the asymptotic growth of $C$ is dictated by $c_2$, which behaves asymptotically like $|x_1|^{2\delta_2'}$,  the uniqueness properties  are independent of the magnitude of $\delta_2'$.

\begin{prop}\label{pubd4.2}$\;$Let $H$ denote the Gru\v{s}in operator defined above.
\begin{tabel}
\item
If $\delta_1\in[0,1/2\rangle$ then $H$ is not Markov unique and consequently not $L_1$-unique.
\item
If $\delta_1\in[1/2,1\rangle$ then $H$ is $L_1$-unique and consequently Markov unique.
\end{tabel}
\end{prop}
\proof\
The first statement of the proposition follows from the observations at the end of Section~6 in \cite{RSi2} 
and in particular from Proposition~6.10.

Some care has to be taken in comparing the current statements
with those of \cite{RSi2}.
The operator $H$ is defined on $C_c^\infty(\Ri^2\backslash\{x_1=0\})$ but the operator $H_\delta$  studied in \cite{RSi2} corresponds to the extension of $H$ to $C_c^\infty(\Ri^2)$.
The Friedrichs' extension $H_{\delta,D}$ of  $H_\delta$ is the self-adjoint extension $H_N$  of $H$ which satisfies the  Neumann-type boundary condition $(c_1\partial_1)(0_+, x_2)=(c_1\partial_1)(0_-, x_2)$ on the line of degeneracy $x_1=0$. 
The  Friedrichs' extension $H_D$ of 
 $H$ is, however, the self-adjoint extension with the Dirichlet-type  boundary condition  $\varphi(0_+,x_2)=\varphi(0_-, x_2)$.
 If $\delta_1\in[0,1/2\rangle$ then these extensions are distinct  and, in addition there are extensions with  analogous Robin boundary conditions sandwiched between the minimal extension $H_N (=H_{\delta,D})$ and  the maximal extension $H_D$.
But if $\delta_1\in[1/2,1\rangle$ then $H_N=H_D$ and all the  operators coincide
(see \cite{RSi2}, Proposition~6.10).

The proof of $L_1$-uniqueness for $\delta_1\in[1/2,1\rangle$ is by reasoning similar to that used to prove
 Proposition~\ref{pubc3.1} and it  does not require an upper bound on  $\delta_2'$.
The argument follows the lines of the proof 
of Theorem~6.1 in \cite{RSi2}, details of which are given in \cite{RSi2a}.
First,  Markov uniqueness follows from  Proposition~6.10 of \cite{RSi2}.
Secondly, one deduces  that $H$ is conservative by  the arguments given in \cite{RSi2a}.
The semigroup $S$ generated by $H_N (=H_D)$ is approximated on $L_2(\Omega)$ by semigroups  $S^{(N,\varepsilon)}$ 
generated by the  Gru\v{s}in operators with coefficients $(C\wedge NI)+\varepsilon I$.
 Then $S$ and $S^{N,\varepsilon}$ satisfy $L_2$-off-diagonal bounds with respect to the corresponding Riemannian distances  by \cite{RSi2}, Proposition~4.1.
 But if $N\geq 1\geq \varepsilon>0$ then these distances are all larger than the Riemannian distance $d_1(\,\cdot\,;\,\cdot\,)$ corresponding to the Gru\v{s}in operator with coefficients $(c_1+1, c_2+1)$.
 Therefore $S$ and the approximants $S^{(N,\varepsilon)}$ all satisfy $L_2$-off-diagonal bounds with respect to $d_1(\,\cdot\,;\,\cdot\,)$.
Since the operator with coefficients $(c_1+1, c_2+1)$ has $\delta_1=0=\delta_2$ it follows that 
$d_1(\,\cdot\,;\,\cdot\,)$ is independent of $\delta_1$ and $\delta_2$.

Next let $B_{1,r}=\{x\in \Ri^2:d_1(0\,;x)<r\}$.
Then if $\varphi\in L_2(B_{1,m})$ it follows by $L_2$-off-diagonal bounds, similar to those of Lemma~\ref{lclo4},  as in the proof of Lemma~\ref{lubc5} that 
 \[
|(\one_{(B_M)^{\rm c}}, S^{(N,\varepsilon)}_t\varphi)|\vee |(\one_{(B_M)^{\rm c}}, S_t\varphi)|\leq 
\sum_{p\geq M}|B_{1,p+1}|^{1/2}
 e^{-d_1(C_p ;B_{1,m})^2 (4t)^{-1}}\|\varphi\|_2
 \]
 where $C_p=B_{1,p+1}\backslash B_{1,p}$.
 But $d_1(C_p ;B_{1,m})\geq p-m$ by the triangle inequality.
 Moreover, it follows from Proposition~5.1 of \cite{RSi2} that there is an $a>0$ such that $|B_{1,p}|\leq a^2\, p^{D'}$
 with $D'=1+(1+\delta_2'-\delta_1')(1-\delta_1')^{-1}$.
 Therefore 
\[
|(\one_{(B_M)^{\rm c}}, S^{(N,\varepsilon)}_t\varphi)|\vee |(\one_{(B_M)^{\rm c}}, S_t\varphi)|\leq 
a\sum_{p\geq M}p^{D'/2}
 e^{-(p-m)^2 (4t)^{-1}}\|\varphi\|_2
 \]
and the estimate is uniform for all  $N\geq1 $ and $\varepsilon\in \langle0,1]$.
Finally the  $S^{(N,\varepsilon)}$ are conservative, since their generators are strongly elliptic,
and they are $L_2$-convergent to $S$.
But 
\[
\|(S^{(N,\varepsilon)}_t-S_t)\varphi\|_1\leq |B_{1,M}|^{1/2}\|(S^{(N,\varepsilon)}_t-S_t)\varphi\|_2+
(\one_{(B_M)^{\rm c}}, S^{(N,\varepsilon)}_t\varphi)+(\one_{(B_M)^{\rm c}}, S_t\varphi)
\]
for all positive $\varphi\in L_2(B_{1,m})$.
Therefore taking the limits $N\to\infty$, $\varepsilon\to0$ and $M\to\infty$ one deduces that the $S^{(N,\varepsilon)}$
are $L_1$-convergent to $S$.
Hence  $S$ is conservative and $H$ is $L_1$-unique by Theorem~\ref{tsm1.1}.
\hfill$\Box$

\vspace{-4mm}

\subsection{Capacity estimates}\label{S4.2}

In this subsection we suppose that $\Omega$ is a strict subset of $\Ri^d$ and examine  the capacity condition
 $\capp_\Omega(\partial\Omega)=0$.
 It follows from the general properties of the capacity  \cite{BH}  or \cite{FOT} that 
  $\capp_\Omega(\partial\Omega)=0$ if and only if  $\capp_\Omega(A)=0$ for all bounded measurable subsets 
  $A$ of $\partial\Omega$.
  Moreover, $\capp_\Omega(A)=0$ implies $|A|=0$.
   The capacity $\capp_\Omega(A)$ depends on two gross features of $A$ and $H_\Omega$,
the dimension of the set and  the order of degeneracy of $H_\Omega$ at~$A$.

First  let  $A\subset \Omega$ be a bounded measurable subset with $|A|=0$ and let $\dim_H(A) $
and $\dim_M(A)$  denote the Hausdorff and  Minkowski dimensions of the set, respectively.
If $A$ is a general measurable subset the corresponding dimensions are defined by $\dim_H(A)=\sup\{\dim_H(B): B\subset A\}$ and $\dim_M(A)=\sup\{\dim_M(B): B\subset A\}$ where the suprema are over all bounded measurable subsets $B$.
In general $\dim_H(A)\leq \dim_M(A)$.

Secondly, the order of degeneracy of the operator $H_\Omega\in \ce_\Omega$ on the bounded set $A$ is defined to be 
the largest  $\gamma_\Omega(A)\geq0$ for which there is an open subset $U$  containing $A$
and an $a>0$  such that $0< C(x)\leq a\,d(x\,;A)^{\gamma_\Omega(A)}I$ for all $x\in U\cap\Omega$ where $d(x\,;A)$ denotes the Euclidean distance of $x$ from $A$.
Again the order of  degeneracy of an unbounded set is defined as a supremum over bounded subsets.

It follows from the proof of Proposition~4.3 in \cite{RSi4} that one has the following property.

\begin{lemma}\label{lsm4.1}
Let  $A\subset \overline\Omega$ be a bounded measurable subset with $|A|=0$.
Assume $H_\Omega\in \ce_\Omega$ is  degenerate of order $\gamma_\Omega(A)$ on $A$.

If $d-2\geq \dim_M(A)-\gamma_\Omega(A)$ then $\capp_\Omega(A)=0$.
In particular if $\gamma_\Omega(A)\geq 2$ then $\capp_\Omega(A)$.
\end{lemma}

The proof of the first statement is similar to the proof of Proposition~4.3 in \cite{RSi4} because the estimates
in the proof of the latter proposition  are all local.
Any possible growth of the coefficients plays no role.
The second statement follows immediately since $\dim_M(A)<d$.

\smallskip

The lemma is not optimal and we next establish a similar but stronger statement involving the Hausdorff dimension.

\begin{prop}\label{pubc4.3}
Let  $A\subset \overline\Omega$ be a bounded measurable subset with $|A|=0$.
Assume $H_\Omega\in \ce_\Omega$ is  degenerate of order $\gamma_\Omega(A)$ on $A$.

If $d-2\geq \dim_H(A)-\gamma_\Omega(A)$ then $\capp_\Omega(A)=0$.
In particular if $d\geq2$ and $\dim_H(A)\leq d-2$ then $\capp_\Omega(A)=0$.
\end{prop}
\proof\
First, if $\gamma_\Omega(A)\geq 2$ then $\capp_\Omega(A)$ by Lemma~\ref{lsm4.1}.
Therefore we assume $\gamma_\Omega(A)<2$.

Secondly, let $B_r=B(y\,;r)$ be a Euclidean ball of radius $r$ centred at $y$ and assume $A\cap B_r\neq\emptyset$.
Then there is an $a>0$ such that $d(x\,;A)\leq a\,r$ for all $x\in B_{2r}$ uniformly in~$y$.
Therefore $\|C(x)\|\leq b\,r^{\gamma(A)}$ for all $x\in B_{2r}\cap \Omega$ with $b=a^{\gamma(A)}$ independent of $y$.
Next fix $\eta_r\in C_c^\infty(B_{2r})$ with $0\leq \eta_r\leq 1$ and $\eta_r=1$ on $B_r$.
One may assume $|\nabla\eta_r|\leq 2\,r^{-1}$ on $B_{2r}\backslash B_r$.
Then
\begin{eqnarray*}
\capp_{\Omega}(B_r\cap\Omega)
\leq\int_{B_{2r}\cap \Omega}\Big(\|C(x)\|\,|(\nabla\eta_r)(x)|^2+|\eta_r|^2\Big)
\leq |B_{2r}|\,(4\,b\,r^{\gamma(A)-2}+1)
\;
\end{eqnarray*}
Since $\gamma(A)<2$ it follows that there is a $c>0$ such that 
$\capp_{\Omega}(B_r\cap\Omega)\leq c\,r^{d+\gamma(A)-2}$
for all $r\leq 1$.
Again the estimate is   uniform in $y$,  the centre point of $B_r$.

Thirdly, let  $B_{r_i}=B(y_i\,;r_i)$ be a countable family of balls with $r_i\leq \delta\leq 1$ such that  $B_{r_i}\bigcap A\neq\emptyset$ and 
$A\subset \bigcup_iB_{r_i}$.
Then
\[
\capp_\Omega(A)\leq \sum_i\capp_\Omega(B_{r_i}\cap\Omega)\leq c\sum_ir_i^{\gamma(A)+d-2}
\]
by the foregoing estimate.
Therefore 
\[
\capp_\Omega(A)\leq c\,H_{\gamma(A)+d-2}(A)
\]
where $H_s$ is the Hausdorff measure.
It follows immediately that $\capp_\Omega(A)=0$ if $\gamma(A)+d-2\geq \dim_H(A)$
or, equivalently, if $d-2\geq \dim_H(A)-\gamma(A)$.
Finally since $\gamma(A)\geq0$ it follows that $\capp_\Omega(A)=0$ whenever $d\geq2$ and $\dim_H(A)\leq d-2$.
\hfill$\Box$

\bigskip

Combination of Theorem~\ref{tsm1.1} and Lemma~\ref{lsm4.1} immediately gives the following criterion for $L_1$-uniqueness.

\begin{prop}\label{pubc4.2}
If the growth condition $(\ref{eubc1.2})$  is satisfied and if 
$d-2\geq \dim_H(A)-\gamma_\Omega(A)$ for each bounded measurable subset $A\subset\partial\Omega$  with $|A|=0$ 
 then $H_\Omega$ is $L_1$-unique.
\end{prop}

Lemma~\ref{lsm4.1}  and Proposition~\ref{pubc4.3} give a quantative assessment of  the two distinct effects which lead to zero capacity.
First  if the operator $H_\Omega$ has a second-order degeneracy on the set $A$, i.e.\ if $\gamma_\Omega(A)\geq 2$, then  $\capp_\Omega(A)=0$ independently of the dimension of $A$.
Secondly, if  $\dim_H(A)\leq d-2$ then $\capp_\Omega(A)=0$
independently of the order of degeneracy of $H_\Omega$ on $A$.
Alternatively, if $\dim_H(A)\leq d-1$
then a first-order degeneracy,  $\gamma_\Omega(A)\geq 1$,  is sufficient to ensure that  $\capp_\Omega(A)=0$.
If, for example,  $A$ is Lipschitz continuous then $\dim_H(A)=d-1$ and this last case is applicable.

\smallskip

In the first of the foregoing cases one even has a simple criterion for $L_2$-uniqueness.
 
\begin{lemma}\label{cubc4.2}
Assume $|\partial\Omega|=0$.
Let  $H_\Omega\in\ce_\Omega$.  
If $\gamma_\Omega(A)\geq 2$ for all bounded measurable $A\subset\partial\Omega$ and
 $\mu(s)\to\infty$ as $s\to\infty$ then
$H$ is $L_2$-unique.
\end{lemma}
\proof\
Let $\rho_n$ denote the functions introduced in the proof of \ref{tsm1.1-4}$\Rightarrow$\ref{tsm1.1-3} in Theorem~\ref{tsm1.1}.
Since $\mu(s)\to\infty$ as $s\to\infty$ it follows that $\rho_n$ converges pointwise to $\one_\Omega$ as $n\to\infty$.
Moreover, $\|\Gamma(\rho_n)\|_\infty\leq b^2\,n^{-2}$.
Next define $\chi_n$ on $[0,\infty\rangle$ by $\chi_n(s)=1$ if $s\in[0,n^{-1}\rangle$, $\chi_n(s)=-\log s/\log n$ if $s\in[n^{-1},1]$ and $\chi_n(s)=0$ if $s\geq1$.
Then define $\xi_n$ on $\Omega$ by $\xi_n(x)=\chi_n(d(x\,;\partial\Omega))$. 
Finally define $\eta_n$ by $\eta_n=\rho_n\,(\one_\Omega-\xi_n)$.
It follows that $\eta_n$ converges pointwise to $\one_\Omega$ as $n\to\infty$.
In addition
\[
\|\Gamma(\eta_n)\|_\infty\leq 2\,\|\Gamma(\rho_n)\|_\infty+2\,a\,\|\Gamma(\xi_n)\|_\infty\leq 
2\,b^2\,n^{-2}+2\,a\,(\log n)^{-2}\to 0
\]
as $n\to\infty$
where we have used $\|\Gamma(\xi_n)\|_\infty\leq a\,(\log n)^{-2}$.
The latter estimate follows from the degeneracy assumption.
Hence $H_\Omega$ is $L_2$-unique by Proposition~6.1 of \cite{RSi4}, with $p=2$, and a standard regularization argument.
\hfill$\Box$

\bigskip

Note that in Lemma~\ref{cubc4.2} there is no restraint on the dimension of the boundary $\partial\Omega$.
Therefore the conclusion is valid for sets $\Omega$ with arbitrarily rough boundaries, in particular for fractal boundaries.
The result is, however, not optimal if the boundary is smooth.
Indeed if $d=1$ and $\Omega=\langle0,\infty\rangle$ then 
a degeneracy of order $1$ at the origin is necessary and sufficient for $L_1$-uniqueness and a degeneracy 
of order $3/2$ is necessary and sufficient for $L_2$-uniqueness (see, \cite{RSi3}, Theorem~2.7).

\vspace{-4mm}

\subsection{Negligible sets}\label{S4.3}

The foregoing discussion indicates that sets of Hausdorff dimension lower than $d-2$ are insignificant for 
$L_1$-uniqueness.
In this subsection we verify that this is indeed the case  for  non-degenerate operators and also establish that sets with Hausdorff dimension  less than $d-4$ are negligible for $L_2$-uniqueness.

\begin{prop}\label{pubc4.5} 
Assume $H\in \ce_{\Ri^d}$ with  $d\geq2$ and let $\Gamma\subset \Ri^d$ be a closed subset with $ |\Gamma|=0$.
Further assume the growth condition $\int_0^\infty ds\, s^{d/2}\, e^{-\lambda\,\mu(s)^2}<\infty $ for one  $\lambda>0$.

The following conditions are equivalent:
\begin{tabel}
\item\label{pubc4.5-1}
$\dim_H(\Gamma)\leq d-2$,
\item\label{pubc4.5-2}
$C_c^\infty(\Ri^d\backslash\Gamma)$ is an  $L_1$-core of $H$.
\end{tabel}
\end{prop}
\proof\
The proof  is based on the observation that both conditions of the proposition are equivalent to the capacity
of the set $\Gamma$ being zero. 
But there are three different capacities involved in the argument.

Let $A\subset \Gamma$ be a bounded measurable subset.
Then we define $\capp(A)$ as the capacity of the set measured with respect to $H$.
Explicitly
\begin{eqnarray}
\capp(A)=\inf\Big\{\;\|\psi\|_{D(h)}^2&&\;: \;\psi\in D(h) \mbox{ and  there exists   an open set  }\nonumber\\[-5pt]
&& U\subset \Ri^d
\mbox{ such that } U\supseteq A
\mbox{ and } \psi\geq1 \mbox{ a.\ e.\ on } U \;\Big\}
\label{eubc4.1}
\end{eqnarray}
with $h$ the closed quadratic form associated with $H$.
Next let $\Omega=\Ri^d\backslash \Gamma$.
Set  $\cd=C_c^\infty(\Ri^d\backslash\Gamma)$ and 
 $H_\Omega=H|_\cd$. 
 Then  $H_\Omega\in \ce_\Omega$.
Let  $\capp_\Omega(A)$ denote the capacity measured with respect to $H_\Omega$.
Thus $\capp_\Omega(A)$ is given by (\ref{eubc4.1}) with $h$ replaced by $h_{\Omega, N}$.
Since $h_{\Omega,N}\supseteq h$ it follows that $\capp_\Omega(A)\leq \capp(A)$.
But both capacities can be calculated with functions which are equal to one in an open neighbourhood
of $A$ and on such functions the two forms coincide.
Therefore $\capp_\Omega(A)= \capp(A)$.

Next define $\capp_{1,2}(A)$ as the capacity of the set $A$ measured with respect to the Laplacian.
This latter capacity is defined by (\ref{eubc4.1}) but  $D(h)$  is replaced by  $W^{1,2}(\Ri^d)$.
Now assume  $\capp(A)=0$.
Then $ C_c^\infty(\Ri^d)$ is a core of $h$, by definition.
Therefore   there exist a sequence $\chi_n\in C_c^\infty(\Ri^d)$ and a decreasing sequence of bounded open subsets $U_n\supset A$ such that $0\leq \chi_n\leq 1$, $\chi_n=1$ on $U_n$ and  $h(\chi_n)+\|\chi_n\|_2^2\leq n^{-1}$.
Then fix an $\eta\in C_c^\infty(\Ri^d)$ such that $0\leq \eta\leq 1$ and $\eta=1$ on $U_1$ and hence on each $U_n$.
Set  $\varphi_n=\chi_n\eta$.
It follows that  $\varphi_n\in D(h)$, $0\leq \varphi_n\leq 1$, $\varphi_n=1$ on $U_n$ and 
$h(\varphi_n)+\|\varphi_n\|_2^2\leq a\,n^{-1}$ with $a=2\,(\|\nabla\eta\|_\infty^2+1)$.
Moreover, $\supp\varphi_n\subseteq K$ for all $n$.
But it follows from strict positivity of the matrix of coefficients $C$ that there exists a $ \mu_K>0$ such that 
\[
 \|\varphi\|_{D(h)}^2\geq \mu_K\,\|\varphi\|^2_{W^{1,2}(\Ri^d)}
\]
for all $\varphi\in W^{1,2}(K)$.
Therefore $\capp_{1,2}(A)=0$.

\smallskip
After these preliminaries we turn to the proof of equivalence of the conditions of the proposition.
We prove that both conditions are equivalent to $\capp(\Gamma)=0$.

\smallskip
\noindent\ref{pubc4.5-1}$\;\Leftrightarrow\capp(\Gamma)=0.\;$
It suffices to prove that Condition~\ref{pubc4.5-1} is equivalent to $\capp(A)=0$ for all bounded measurable subsets $A$ of $\Gamma$.
But Proposition~\ref{pubc4.3}   establishes that Condition~\ref{pubc4.5-1} implies $\capp_\Omega(A)=0$ or, equivalently, 
$\capp(A)=0$.
Conversely, $\capp(A)=0$ implies  $\capp_{1,2}(A)=0$ by the foregoing discussion.
But  $\capp_{1,2}(A)=0$ implies $\dim_H(\Gamma)\leq d-2$ by standard properties of the Laplacian (see, for example,  Theorem~2.26 of \cite{HKM}, Corollary~5.1.15 of \cite{AdH} or Section~2.1.7 of \cite{MZ1}).

\smallskip

\noindent\ref{pubc4.5-2}$\;\Leftrightarrow\capp(\Gamma)=0.\;$
First,  
observe that $H$ is  $L_1$-unique  by the second statement of Proposition~\ref{pubc4.1}.
But $H$ is $L_1$-unique if and only if   $\overline H^{\scriptscriptstyle 1}$ is the generator of an $L_1$-continuous semigroup.
Now suppose Condition~\ref{pubc4.5-2} is valid.
Then  $\overline H^{\scriptscriptstyle 1}=\overline{H_\Omega}^{\scriptscriptstyle 1}$.
Therefore $\overline{H_\Omega}^{\scriptscriptstyle 1}$ is the generator of an $L_1$-continuous semigroup.
Consequently,  $H_\Omega$ is  $L_1$-unique.
But  $L_1$-uniqueness of $H_\Omega$ is equivalent to  $\capp_\Omega(\Gamma)=0$ by Theorem~\ref{tsm1.1} which in turn is equivalent to  $\capp(\Gamma)=0$.

Conversely, if $H_\Omega$ is $L_1$-unique then  $\overline{H_\Omega}^{\scriptscriptstyle 1}$  is the generator of an $L_1$-continuous semigroup.
But  $\overline H^{\scriptscriptstyle 1}$ is also a generator and 
$\overline H^{\scriptscriptstyle 1}\supseteq\overline{H_\Omega}^{\scriptscriptstyle 1}$
by definition.
Since a semigroup generator cannot have a proper generator extension one must have
$\overline H^{\scriptscriptstyle 1}=\overline{H_\Omega}^{\scriptscriptstyle 1}=\overline {H|_\cd}^{\scriptscriptstyle 1}$.
Thus $\cd$ is an $L_1$-core of $H$.
\hfill$\Box$

\bigskip

It follows from the assumptions of Proposition~\ref{pubc4.5} and the first statement of Proposition~\ref{pubc4.1}
that $H$ is $L_2$-unique, i.e.\ $H$ is essentially self-adjoint.
Then the closed  form $h$ associated with $H$ is the form of the $L_2$-closure 
$\overline H$ of $H$.
Therefore $h(\varphi)=\|{\overline H}^{1/2}\varphi\|_2^2$ for all $\varphi\in D( h)=D({\overline H}^{1/2})$.
This observation provides a relation between $L_1$- and $L_2$-cores.

\begin{cor}\label{cubc4.1}
Adopt the assumptions of Proposition~$\ref{pubc4.5}$ .

The following conditions are equivalent:
\begin{tabel}
\item\label{cubc4.1-1}
$C_c^\infty(\Ri^d\backslash\Gamma)$ is an  $L_1$-core of $H$,
\item\label{cubc4.1-2}
$C_c^\infty(\Ri^d\backslash\Gamma)$ is an  $L_2$-core of ${\overline H}^{1/2}$.
\end{tabel}
\end{cor}
\proof\
Again set  $\Omega=\Ri^d\backslash \Gamma$, $\cd=C_c^\infty(\Ri^d\backslash\Gamma)$
and $H_\Omega=H|_\cd$.
Then $H_\Omega\in \ce_\Omega$.
Moreover, it follows from the proof of Proposition~\ref{pubc4.5} that Condition~\ref{cubc4.1-1} is equivalent to 
the condition $\capp_\Omega(\Gamma)=0$.
But Condition~\ref{cubc4.1-2} is also equivalent to 
this capacity condition  by the proof of Theorem~\ref{tsm1.1}.
\hfill$\Box$

\bigskip

There is also an $L_2$-version of Proposition~\ref{pubc4.5}.

\begin{prop}\label{pubc4.6} 
Assume $H\in \ce_{\Ri^d}$ with  $d\geq4$ and let $\Gamma\subset \Ri^d$ be a closed subset with $ |\Gamma|=0$.
Further assume that $\mu(s)\to\infty $ as $s\to\infty$.

The following conditions are equivalent:
\begin{tabel}
\item\label{pubc4.6-1}
$\dim_H(\Gamma)\leq d-4$,
\item\label{pubc4.6-2}
$C_c^\infty(\Ri^d\backslash\Gamma)$ is an  $L_2$-core of $H$.
\end{tabel}
\end{prop}
\proof\
It follows from the first statement of Proposition~\ref{pubc4.1} that $H$ is essentially self-adjoint.
Then if $A\subset \Ri^d$  is a measurable subset we define  the capacity $\Capp(A)$ associated with the 
self-adjoint $L_2$-closure $\overline H$ by
\begin{eqnarray}
\Capp(A)=\inf\Big\{\;\|\psi\|_{D(\overline H)}^2\;: \;&&\psi\in D(\overline H) \mbox{ and  there exists   an open set  }\nonumber\\[-5pt]
&& U
\mbox{ such that } U\supseteq A
\mbox{ and } \psi\geq1 \mbox{ a.\ e.\ on } U\;\Big\}
\;.\label{eubc4.11}
\end{eqnarray}
We will argue that both conditions of the proposition are equivalent to $\Capp(\Gamma)=0$ or, equivalently, $\Capp(A)=0$
for all bounded measurable subsets $A$ of $\Gamma$.

If  $A$ is bounded and $\Capp(A)=0$ then 
there exist a sequence $\chi_n\in C_c^\infty(\Ri^d))$ and a decreasing sequence of bounded open subsets $U_n\supset A$ such that $0\leq \chi_n\leq 1$, $\chi_n=1$ on $U_n$ and $\|\chi_n\|_{D(\overline H)}^2\leq n^{-1}$.
But the sequence $\chi_n$ can be modified by a variation of the argument used in the proof of  Proposition~\ref{pubc4.5} to yield a sequence with similar properties but with each element of the sequence supported by a fixed compact set.
Explicitly, fix  $\eta\in C_c^\infty(\Ri^d)$ 
such that  $0\leq \eta\leq 1$ and $\eta=1$ on $U_1$ and hence on each $U_n$.
Then set  $\varphi_n=\chi_n\eta$.
It follows that $0\leq \varphi_n\leq 1$, $\varphi_n=1$ on $U_n$  and $\|\varphi_n\|_2\leq \|\chi_n\|_2$.
But 
\[
\overline H\varphi_n= (\overline H\chi_n)\eta+\chi_n(\overline H\eta)+2\,\Gamma(\chi_n,\eta)
\]
where $\Gamma(\,\cdot\,;\,\cdot\,)$ is the {\it carr\'e du champ} associated with $H$.
Therefore 
\[
\|\overline H\varphi_n\|_2^2\leq 3\,(\|\overline H\eta\|_\infty^2+\|\eta\|_\infty^2)\,\|\chi_n\|_{D(\overline H)}^2
+ 12\int \Gamma(\chi_n;\eta)^2
\;.
\]
But
\[
\int \Gamma(\chi_n;\eta)^2\leq \int \Gamma(\chi_n)\,\Gamma(\eta)
\leq |\Gamma(\eta)|_\infty\, h(\chi_n)\leq 2^{-1}|\Gamma(\eta)|_\infty \|\chi_n\|_{D(\overline H)}^2
\;.
\]
Combining these estimates one deduces that there is an $a>0$ such that 
$\|\varphi_n\|_{D(\overline H)}^2\leq a\,\|\chi_n\|_{D(\overline H)}^2\leq a\,n^{-1}$ for all $n$.

Next  since the coefficients $c_{ij}\in W^{1,\infty}_{\rm loc}(\Omega)$ and  $C=(c_{ij})>0$ it follows that there exists a 
$\mu_K>0$ such that 
\[
 \|\varphi\|_{D(\overline H)}^2\geq \mu_K\,\|\varphi\|^2_{W^{2,2}(\Ri^d)}
\]
for all $\varphi\in C_c^\infty(K)$ (see, for example, the appendix of  \cite{RSi4}).
Therefore replacing $\varphi$ by $\varphi_n$ and taking the limit $n\to\infty$ one deduces that $\capp_{2,2}(A)=0$
where $\capp_{2,2}$ is the capacity measured with respect to the $W^{2,2}(\Ri^d)$-norm, i.e.\
$\capp_{2,2}(A)$ is given by  (\ref{eubc4.11}) but with $D(\overline H)$ replaced by $W^{2,2}(\Ri^d)$.

\smallskip
\noindent\ref{pubc4.6-1}$\;\Leftrightarrow\Capp(\Gamma)=0.\;$
It suffices to prove that Condition~\ref{pubc4.6-1} is equivalent to $\Capp(A)=0$ for all bounded measurable subsets $A$ of $\Gamma$.
But  a slight variation of the proof of Proposition~\ref{pubc4.3}   establishes that Condition~\ref{pubc4.6-1} implies $\Capp(A)=0$.
Indeed define $B_r$ and $\eta_r$ as in the proof of Proposition~\ref{pubc4.3}.
One may also assume that $|\Delta\eta_r|\leq a\,r^{-2}$ on $B_{2r}\backslash B_r$.
Then one estimates straightforwardly that there are  $b,b'>0$ such that
\[
\Capp(B_r)\leq b\int_{B_{2r}}\Big(|\Delta\eta_r|^2+|\eta_r|^2\Big)\leq b'\,(r^{d-4}+r^d)\leq 2\,b'\,r^{d-4}
\]
for all $r\leq1$.
The rest of the proof remains unchanged.

Conversely, if $\Capp(A)=0$ then $\capp_{2,2}(A)=0$ by the foregoing discussion.
But then  $\dim_H(A)\leq d-4$ by  another application of  Theorem~2.26 of \cite{HKM}.

\smallskip

\noindent\ref{pubc4.6-2}$\;\Leftrightarrow\Capp(\Gamma)=0.\;$
First suppose Condition~\ref{pubc4.6-2} is valid. 
Secondly,  fix $\psi\in C_c^\infty(\Ri^d)\subset D(H)$ with $\psi=1$ on an open neighbourhood  of $A\subset \Gamma$.
Then, by \ref{pubc4.6-2},  there is a sequence $\psi_n\in C_c^\infty(\Ri^d\backslash \Gamma)$ such that 
$\|\psi-\psi_n\|_{D(\overline H)}\to0$ as $n\to\infty$.
Set $\varphi_n=\psi-\psi_n$.
It follows that $\varphi_n\in D(H)$, $\varphi_n=1$ on an open neighbourhood $U_n$ of $A$ and 
$\|\varphi_n\|_{D(\overline H)}\to0$ as $n\to\infty$.
Therefore $\Capp(A)=0$.
Since this holds for an arbitrary  bounded subset $A$ of $\Gamma$  it follows that $\Capp(\Gamma)=0$.

Conversely, suppose $\Capp(\Gamma)=0$.
Therefore $\Capp(A)=0$ for each bounded measurable subset $A$ of $\Gamma$.
Then since $C_c^\infty(\Ri^d)$ is a core of $H$, by definition, there exist a sequence $\chi_n\in C_c^\infty(\Ri^d)$ and  a sequence of open subsets $U_n$ of $A$ such that $\|\chi_n\|_{D(\overline H)}\to0$ as $n\to\infty$.
Now fix $\psi\in C_c^\infty(\Ri^d)$ and set $\psi_n=(\one-\chi_n)\psi$.
It follows that  $\psi_n\in C_c^\infty(\Ri^d\backslash \Gamma)$.
Moreover, $\psi-\psi_n=\chi_n\psi$.
But arguing as in the proof of Proposition~\ref{pubc4.5} one has
\[
H(\psi-\psi_n)=(H\chi_n)\psi+\chi_n(H\psi)+2\,\Gamma(\chi_n\,;\psi)
\]
and consequently
\[
\|H(\psi-\psi_n)\|_2^2\leq 3\,\Big((\|H\psi\|_\infty^2+\|\psi\|_\infty^2)+2\,\|\Gamma(\psi)\|_\infty)\Big)\,\|\chi_n\|_{D(\overline H)}^2
\;.
\]
Since $\|\psi-\psi_n\|_2\leq \|\psi\|_\infty\|\chi_n\|_2$ it follows that $\|(\psi-\psi_n)\|_{D(\overline H)}\to0$ as $n\to\infty$.
Therefore each $\psi\in C_c^\infty(\Ri^d)$ can be approximated by a sequence  $\psi_n\in C_c^\infty(\Ri^d\backslash \Gamma)$ in the $D(\overline H)$-graph norm.
But as $C_c^\infty(\Ri^d)$ is a core of $H$ one concludes that $C_c^\infty(\Ri^d\backslash \Gamma)$ is also a core.
\hfill$\Box$

%\bibliography{refbib}

%
%\newpage

\end{document}